\documentclass[preprint]{article}

\usepackage{graphicx}
\usepackage{bm}
\usepackage{amsmath}
\usepackage{amssymb}
\usepackage{color}
\usepackage{epstopdf}
\def\ar#1{\textcolor{black}{#1}}

\newcommand{\bc}{\begin{center}}
\newcommand{\ec}{\end{center}}
\newcommand{\bd}{\begin{displaymath}}
\newcommand{\ed}{\end{displaymath}}
\newcommand{\be}{\begin{equation}}
\newcommand{\ee}{\end{equation}}
\newcommand{\ben}{\begin{eqnarray}}
\newcommand{\een}{\end{eqnarray}}
\newcommand{\benn}{\begin{eqnarray*}}
\newcommand{\eenn}{\end{eqnarray*}}
\newcommand{\bi}{\begin{itemize}}
\newcommand{\ei}{\end{itemize}}

\newcommand{\bfa}[1]{\mbox{\boldmath $ #1 $}}

\newcommand{\vE}{\bfa{E}}

\newcommand{\vl}{{\bfa l}}
\newcommand{\vm}{{\bfa m}}
\newcommand{\vn}{{\bfa n}}

\newcommand{\vQ}{{\bfa Q}}

\bmdefine{\ehat}{e}
\bmdefine{\nhat}{n}
\bmdefine{\bmD}{D}
\bmdefine{\bmE}{E}
\bmdefine{\vlambda}{\lambda}
\bmdefine{\vdelta}{\delta}
\newcommand{\eps}{\varepsilon}
\bmdefine{\bfeps}{\eps}

\title
{A Moving Mesh Method for Modelling Defects in Nematic Liquid Crystals}

\author{Craig S. MacDonald, John A. Mackenzie, Alison Ramage\thanks{
Department of Mathematics and Statistics, University of 
Strathclyde, Glasgow G1 1XH, Scotland.}}

\begin{document}
\maketitle

\begin{abstract}
The properties of liquid crystals can be modelled using an order parameter 
which describes the variability of the local orientation of rod-like molecules. Defects in the 
director field can arise due to external factors such as applied electric or magnetic fields, or 
the constraining geometry of the cell containing the liquid crystal material. Understanding the 
formation and dynamics of defects is important in the design and control of liquid crystal devices,
and poses significant challenges for numerical modelling.
In this paper we consider the numerical solution of a $\bfa{Q}$-tensor model of a nematic liquid 
crystal, where defects arise through rapid changes in the $\bfa{Q}$-tensor over a very small 
physical region in relation to the dimensions of the liquid crystal device. The efficient solution 
of the resulting six coupled partial differential equations is achieved using a finite 
element based adaptive moving mesh approach, where an unstructured triangular mesh 
is adapted towards high activity regions, including those around defects. Spatial convergence 
studies are presented using a stationary defect as a model test case, and the adaptive method is 
shown to be optimally convergent using quadratic triangular finite elements. The full effectiveness 
of the method is then demonstrated using a challenging two-dimensional dynamic Pi-cell 
problem involving the creation, movement, and annihilation of defects.

\end{abstract}

\section{Introduction}

The orientational properties of liquid crystal materials can be manipulated by applying an electric 
or magnetic field, leading to particular characteristics of the reflection and transmission of 
light waves. These effects make liquid crystals key materials in the
construction of a broad range of commonly-used display devices, such as 
the Twisted Nematic Device (TND) \cite{schadt:71}, the Pi-cell \cite{bos:84} and the 
Zenith Bistable-Device (ZBD) \cite{bryanbrown:95,newton:97}. More recently, there has been growing 
interest in liquid crystals in a wider context. Examples include active liquid crystals 
\cite{majumdar:14} (which are relevant to natural applications such as modelling cytoskeletal 
structure in cell biology or animal flocking as well as in synthetic manufacture of colloids and 
granular matter), liquid crystal shells and drops \cite{lopez:11}, and materials design and 
self-assembly of ordered fluids \cite{Tschierske:12}. Because of their importance in these and
other technological applications, there is a great deal of interest in modelling the 
properties of liquid crystals mathematically.

 The most commonly-used continuum models 
utilise one or more unit vector fields as state variables. For the uniaxial nematic phase, which is 
the simplest and most common liquid crystal phase, the orientation of the molecules is represented 
by a unit vector denoting the direction in which their main axis points. This is known as the
\textit{liquid crystal director} and is traditionally denoted by $\vn$. More generally, 
taking $\vn$ and $-\vn$ to be equivalent, the average molecular orientation can be represented 
by an order tensor, usually denoted by $\vQ$. This tensor can be written as
\be
\vQ = S\left(\vn \otimes \vn-\frac{1}{3}\bfa{I}\right) + T(\vm \otimes \vm - \vl\otimes\vl),
\label{Qform}
\ee
where $S$ and $T$ are scalar order parameters, $\{\vl,\vm,\vn\}$ is a set of orthonormal 
directors, and $\bfa{I}$ is the identity (see, for example, \cite{sonnet:12}). Note that 
the uniaxial case can be recovered by setting $T=0$. 

In this paper, we propose an efficient numerical method for computing the orientational 
state of a nematic liquid crystal based on a $\vQ$-tensor model. In particular, we focus on 
tracking the movement of \textit{defects} in the material, that is, local regions (of 
\textit{point}, \textit{line} or 
\textit{wall} type) where the symmetry of the ordered material is broken. 
The switching behaviour of liquid crystal material between two equilibrium states (by means of an 
applied field), which is the basis of most liquid crystal devices, is strongly influenced by 
the existence of such defects, so it is important to be able to model these features accurately.
Our use of the $\vQ$-tensor (as opposed to a director-based) model in this paper is driven by the 
fact that in this formulation, topological defects do not appear as mathematical singularities. 

A $\vQ$-tensor theory of nematic liquid crystals, which allows for changes in the scalar 
order parameters, has been developed from the theory of Landau by de Gennes \cite{degennes:72}. 
Minimisation of the total free energy in the case of a nematic liquid crystal coupled with an 
applied electric field leads to a set of six coupled partial differential equations (PDEs) for the 
five degrees of freedom of the order parameter tensor $\vQ$ and the electric potential $U$, which 
poses a challenge for numerical solvers. Furthermore, additional physical features such as flow and temperature change require the $\vQ$-tensor equations are coupled to the Navier-Stokes and 
energy equations. Even in the absence 
such additional complications, the $\bfa{Q}$-tensor equations are difficult to solve numerically due 
to their highly non-linear nature. Also, the defects mentioned above induce distortion of the 
director over very small length scales as 
compared to the size of the cell. It can therefore be difficult to accurately represent their 
nature and behaviour with a standard numerical model. The large discrepancies in length and time 
scales which occur mean that numerical difficulties are even more acute for models of dynamic 
problems involving the movement of defects, such as the Pi-cell problem studied in \S\ref{Picell}.

For identifying static equilibrium states, relatively straightforward numerical methods are often 
good enough (see, for example, \cite{davis:98, gartland:02, lee:02,mottram:04, 
schopohl:87,sonnet:95}). There have also been several studies using more sophisticated adaptive 
techniques. These include the $h$ (grid parameter) and $p$ (degree of basis function) adaptive 
finite element methods presented in \cite{newton:11, fukuda:01,fukuda:02, james:06}. 
Additional methods have been proposed based on moving meshes 
\cite{abl10,abl11,macdonald:11,macdonald:15,rn07,rn08}, 
which move existing mesh points so as to cluster them in areas of large solution error
whilst maintaining the same mesh connectivity. These techniques are particularly appropriate for 
resolving localised solution singularities such as defects, as maintaining a fixed connectivity is 
very efficient in terms of computing time (as opposed to adding or removing grid points in areas 
of interest). Also, for transient problems, it is sometimes possible to use bigger time steps if 
the solution remains almost stationary relative to the moving mesh frame of reference. 
This motivates our use of adaptive moving mesh
techniques to capture defect structure and track defect movement within the cell. 

In \cite{macdonald:15}, we proposed a robust and efficient numerical scheme for solving the system
of six coupled partial differential equations which arises when using $\vQ$-tensor theory to model the
behaviour of a nematic liquid crystal cell under the influence of an applied electric field in one 
space dimension. The numerical method uses a moving mesh partial differential equation (MMPDE) approach to
generate an adaptive mesh which accurately resolves important solution features. In this paper, we 
extend this adaptive moving mesh strategy to solve liquid crystal problems in two dimensions. 
This involves addressing a number of significant new challenges, including 
the choice of appropriate adaptivity criteria for problems with moving
singularities, the efficient solution of the large systems of highly
non-linear algebraic equations arising after discretisation, and how to deal with the 
creation and annihilation of defects in a realistic model.

The remainder of the paper is structured as follows. In \S\ref{sec:2}, we give a brief overview of 
the derivation of the physical PDEs arising from the $\vQ$-tensor framework coupled with an 
applied electric field, along with some details of their finite element discretisation when an 
adaptive moving mesh is utilised. In \S\ref{moving},  the details of the two-dimensional moving mesh PDE are given.
We consider a number of different mesh adaptivity criteria through the use of monitor functions, and present 
a series of numerical experiments which indicate that monitor functions based on a local measure of 
biaxiality perform well.  We then apply the biaxiality-based monitor function 
to a problem first presented by Bos \cite{bos:07}: a two dimensional Pi-cell 
problem with a sinusoidal perturbation across the centre of the cell. This is a 
dynamic two-dimensional version of the test problem described in \cite[\S1.2]{macdonald:15}.

\section{Derivation and discretisation of physical PDEs}
\label{sec:2}
\subsection{Derivation of physical PDEs}
\label{sec:qtensor}

To characterise the molecular alignment of a nematic liquid crystal, 
we define a \textit{uniaxial} $\vQ$-tensor using a local ensemble average of the molecular axes as
\be
\bfa{Q}=\sqrt{\frac{3}{2}}\left\langle {\bfa u} \otimes {\bfa u}
-\frac{1}{3}\bfa{I}\right\rangle
\label{Qu}
\ee
(see, e.g., \cite[\S2.1.2]{degennes:physics2}). 
The unit vectors ${\bfa u}$ lie along the molecular axes and 
the angle brackets denote the ensemble averaging: 
the factor $\sqrt{3/2}$ is included for convenience so that, for a uniaxial state with director 
$\bfa{n}$ and scalar order parameter $S$, $\mathrm{tr}(\bfa{Q}^2)=S^2$. 
The tensor (\ref{Qu}) has five degrees of freedom and is symmetric and traceless,
so it can be represented in matrix form as
\be
 \bfa{Q}= \left[\begin{array}{ccc}
q_1\quad&q_2\quad&q_3\\q_2\quad&q_4\quad&q_5\\q_3\quad&q_5\quad&-q_1-q_4
\end{array}\right],
\label{Qmat}
\ee
where each of the five quantities $q_i$, $i=1,\ldots,5$, is a function of 
the spatial coordinates and time. Note that the orthonormal eigenvectors of this matrix are 
the vectors $\{\vl, \vm, \vn\}$ used in the representation of the $\vQ$ tensor given in 
(\ref{Qform}).

The globally stable state of a nematic liquid crystal under the 
influence of an applied electric field
corresponds to a minimum point of the free energy. Using Landau-de Gennes theory, in which the free 
energy density is assumed to depend on $\bfa{Q}$ and its gradient, the free energy may be written 
as 
\begin{equation}
 {F}=\int_V (
\mathcal{F}_t(\mathbf{Q})+
\mathcal{F}_e(\mathbf{Q},\nabla \mathbf{Q})+
\mathcal{F}_u(\mathbf{Q},\nabla \mathbf{Q}))\,{\rm 
d}V,
\label{free_en}
\end{equation}
where $\mathcal{F}_t$, $\mathcal{F}_e$, $\mathcal{F}_u$ and $\mathcal{F}_s$
represent the thermotropic, elastic and electrostatic 
terms,
respectively. Note that, as here we only consider problems with fixed (\textit{strong anchoring}) 
boundary conditions, we omit any (constant) surface energy terms.
Expressions for the individual terms in the integrand of (\ref{free_en}) can be derived in a 
variety of different ways: here we expand the thermotropic energy, $\mathcal{F}_t$, up to 
fourth order in $\bfa{Q}$ and the elastic energy, $\mathcal{F}_e$, up to second order in the 
gradient of
$\bfa{Q}$. Details of the resulting expressions can be found in \cite[\S2]{macdonald:15}, along 
with a description of the contribution from the applied electric field, $\mathcal{F}_u$, which 
includes flexoelectricity. As in \cite{macdonald:15}, values for material constants 
throughout this paper are those used in \cite{barberi:04},
which are commensurate with a liquid crystal cell of the 5CB
compound 4-cyano-$4^{\prime}$-n-pentylbiphenyl.

To derive time-dependent PDEs for the quantities $q_i$ in (\ref{Qmat}), we use a
dissipation principle with viscosity coefficient $\nu$ and dissipation function 
\begin{displaymath}
 \mathcal{D}=\frac{\nu}{2}\textrm{tr}\left[\left(\frac{\partial  \textit{\textbf{Q}}}
{\partial
t}\right)^2\right]=\nu(\dot{q}_1\dot{q}_4+\dot{q}_1^2+\dot{q}_2^2+\dot{q}_3^2+\dot
{q}_4^2+\dot{q}_5^2), 
\end{displaymath}
where the dot represents differentiation with respect to time (see, e.g., \cite[eq.\ 
(4.23)]{sonnet:12}). 
For a physical domain with spatial coordinates $\{x_1, x_2, x_3\}$,
this produces a system of equations which can be written as
\begin{equation}
 \frac{\partial
\mathcal{D}}{\partial{\dot{q}}_i}=\nabla\cdot{\hat{\bfa{\Gamma}}}_i-{\hat f}_i
\qquad i=1,\ldots,5,
\label{diss}
\end{equation}
involving the bulk energy $\mathcal{F}_b=\mathcal{F}_t+\mathcal{F}_e+\mathcal{F}_u$, 
where the vector ${\hat{\bfa{\Gamma}}}_i$ has entries
\begin{displaymath}
 ({\hat{\bfa{\Gamma}}}_i)_j=\frac{\partial {\mathcal F}_b}{\partial q_{i,j}},\qquad
q_{i,j}=\frac{\partial q_i}{\partial x_j},\qquad
j=1,2,3,
\end{displaymath}
and  ${\hat f}_i$ is given by
\begin{displaymath}
{\hat f}_i=\frac{\partial {\mathcal F}_b}{\partial q_i}.
\end{displaymath}
To add the coupling with an electric field, $\vE$ say, 
we introduce an additional unknown in the form of a scalar 
electric potential $U$ such that $\vE=-\nabla U$. Assuming that 
there are no free charges, the electric field within the cell can then be found by solving the 
Maxwell equation
\begin{equation}
\nabla\cdot \mathbf{D}=0
\label{eq:max}
\end{equation}
where the specific form of the electric displacement $\mathbf{D}$ can be found in \cite[eqn 
(2.5)]{macdonald:15}.
Minimisation of the total free energy (\ref{free_en}) therefore leads to a set of six coupled 
non-linear PDEs
for the five degrees of freedom of $\bfa{Q}$ and the electric potential $U$.  Specifically,
combining (\ref{diss}) and (\ref{eq:max}) and using some algebraic manipulation for notational 
convenience, we obtain the equations
\begin{subequations}
\begin{equation}
\displaystyle{\frac{\partial
q_{i}}{\partial t}}  =  \nabla\cdot \bfa{\Gamma}_{i}-
f_{i},\;\; \quad i=1,\ldots,5,\label{eq:1Deqs}
\end{equation}
\begin{equation}
\nabla \cdot {\bfa D}=0,
\end{equation}
\label{eq:conduf}
\end{subequations}
where 
\[
{\bfa \Gamma}_{1}=\frac{1}{3\nu}(2{\hat{\bfa \Gamma}}_{1}-{\hat{\bfa \Gamma}}_{4}), \quad \quad
f_{1}=\frac{1}{3\nu}(2{\hat f}_{1}-{\hat f}_{4}),
\]
\[
{\bfa \Gamma}_{4}=\frac{1}{3\nu}(2{\hat{\bfa \Gamma}}_{4}-{\hat{\bfa \Gamma}}_{1}), \quad \quad
f_{4}=\frac{1}{3\nu}(2{\hat f}_{4}-{\hat f}_{1})
\]
and 
\[
{\bfa \Gamma}_{i}=\frac{1}{2\nu}{\hat{\bfa \Gamma}}_{i}, \quad \quad f_{i}=\frac{1}{2\nu}{\hat{f}}_{i}, 
\quad i=2,3,5.
\]
Note that we non-dimensionalise the equations (\ref{eq:conduf}) for computational purposes: details 
of the precise scalings applied in terms of lengths and energies are given in 
\cite[\S2]{macdonald:15}. 

The governing physical PDEs in (\ref{eq:conduf}) now have to be adapted to account 
for the movement of the finite element mesh. To do this, we introduce a family of bijective 
mappings 
\begin{equation}
{\cal A}_{t}: \Omega_{c}\subset \mathbb{R}^2 \rightarrow
\Omega\subset\mathbb{R}^2, \quad \quad {\bfa x}({\bfa \xi},t)={\cal A}_{t}({\bfa \xi}),
\label{eq:map1}
\end{equation}
such that, at a given time $t$ in time domain $T\subseteq\mathbb{R}^{+}$,  the point ${\bfa 
\xi}=(\xi,\eta)$ of a
two-dimensional \ar{computational reference} domain $\Omega_{c}$ is mapped to the point
${\bfa x}=(x,y)$ of the original physical domain $\Omega$. 
The temporal derivative of a mapping $g: \Omega \rightarrow \mathbb{R}$ (from the physical
domain) in the computational frame can then be defined as
\begin{displaymath}
\left .\frac{\partial g}{\partial t}\right |_{ {\bfa \xi}} : \Omega \rightarrow
\mathbb{R}, 
\quad \left .\frac{\partial g}{\partial
t}\right |_{{\bfa \xi}} ({\bfa x},t)= \frac{\partial {\hat
g}}{\partial t}({\bfa \xi},t), \quad {\bfa \xi}={\cal A}^{-1}_{t}({\bfa x}),
\end{displaymath}
where ${\hat g}:\Omega_{c}\times T\rightarrow \mathbb{R}$ is the
corresponding function in the computational domain, that is, ${\hat
g}({\bfa \xi},t)={\hat g}(({\bfa x},t),t)=g({\cal A}_{t}({\bfa \xi}))$.
Defining the \textit{mesh velocity} $\dot{\bfa x}$ as
\begin{displaymath}
\dot{\bfa x}({\bfa x},t)=
\left . \frac{\partial {\bfa x}}{\partial t}
\right |_{{\bfa \xi}}({\cal A}^{-1}_{t}({\bfa x})),
\end{displaymath}
and applying the chain rule for differentiation (with appropriate smoothness assumptions on $g$) 
gives
\begin{displaymath}
\left . \frac{\partial q}{\partial t}\right |_{{\bfa \xi}}= \left .
\frac{\partial q}{\partial t}\right |_{\bfa x} + \dot{{\bfa x}}
\cdot\nabla q,
\label{eq:udot}
\end{displaymath}
which includes an additional convection-like term due to the mesh movement.
Incorporating these changes into (\ref{eq:1Deqs}) gives the final set of six coupled PDEs
\begin{subequations}
\begin{equation}
\left . 
\frac{\partial q_{i}}{\partial t}\right |_{{\bfa \xi}}  -
\dot{\bfa x}\; \cdot\nabla q  =  
\nabla \cdot\Gamma_{i} -f_i
\;\; \quad i=1,\ldots,5,
\label{eq:aleconduf}
\end{equation}
\begin{equation}
\nabla \cdot {\bfa D}=0.
\end{equation}
\label{eq:ale}
\end{subequations}

\subsection{Finite element discretisation}
\label{femethod} 
With a space of time-independent finite element test functions $\hat{v}\in H_{0}^{1}(\Omega_{c})$, 
mesh mapping (\ref{eq:map1}) defines the test space
\[
{\cal H}_{0}(\Omega)=\left \{ v: \Omega \rightarrow \mathbb{R}:
v=\hat{v} \circ {\cal A}_{t}^{-1}, \; \hat{v}\in H_{0}^{1}(\Omega_{c})\right \}.
\]
We denote the approximation spaces with essential boundary conditions on 
$q_{i}$ and $U$ by ${\cal H}_{E_{q}}$ and ${\cal H}_{E_{U}}$, respectively.
In an analogous way to the one-dimensional case described in \cite[\S3.1]{macdonald:15}, 
Reynolds's transport formula can be used to derive the following conservative 
weak formulation of (\ref{eq:ale}):
find $q_{i}\in {\cal H}_{E_{q}}(\Omega)$,  $i=1,\ldots,5$, and 
$U\in {\cal H}_{E_{U}}(\Omega)$ such that $\forall v \in {\cal H}_{0}(\Omega)$ 
\begin{subequations}
\begin{equation}
\frac{d}{dt} \int_{\Omega} q_{i}v  \; {\rm d}{{\bfa x}}-
\int_{\Omega} (\nabla \cdot ({\dot {\bfa x}}q_{i}))\,v\;{\rm d}{{\bfa x}} 
= \int_{\Omega} \Gamma_{i}\cdot \nabla v \; {\rm d}{{\bfa x}}- \int_{\Omega} f_{i} v \; {\rm d}{{\bfa x}}, 
\quad  
\label{eq:wcon}
\end{equation}
\begin{equation}
\int_{\Omega} {\bfa D}\cdot \nabla v \; {\rm d}{{\bfa x}} =  0.
\label{eq:ddismax}
\end{equation}
\label{eq:weak}
\end{subequations}

To discretise (\ref{eq:weak}), we assume that the reference domain $\Omega_{c}$ is covered by a 
fixed triangulation ${\cal T}_{h,c}$ with straight edges, so that $\Omega_{c} = \cup_{K\in {\cal 
T}_{h,c}}K$,  and
introduce the Lagrangian finite element spaces
\begin{eqnarray*}
{\cal L}^{k}(\Omega_{c})&  = & \{ {\hat v}_{h} \in
H^{1}(\Omega_{c})\, : \, {\hat v}_{h}|_{K} \in
\mathbb{P}_{k}(K), \;\; \forall \, K \in {\cal T}_{h,c}  \} \nonumber \\
{\cal L}^{k}_{0}(\Omega_{c})&  = & \{ {\hat v}_{h} \in
H^{1}(\Omega_{c})\, : \, {\hat v}_{h}|_{K} \in {\cal
L}^{k}(\Omega_{c})\, : \, {\hat v}_{h} = 0, \;\; {\bfa \xi}\in
\partial \Omega_{c} \},
\end{eqnarray*}
where ${\mathbb P}_{k}(K)$ is the space of polynomials on $K$ of degree
less than or equal to $k$.   
Using a piecewise linear discretisation of the mesh mapping (\ref{eq:map1}) 
to produce a discrete mapping ${\cal A}_{h,t} \in {\cal L}^{1}(\Omega_{c})$,
finite element spaces on the physical domain $\Omega$ can be defined as
 \begin{eqnarray*}
{\cal L}^{k}(\Omega)&  = & \{
v_{h} \, : \, \Omega\rightarrow {\mathbb R}\,:\, v_{h}= {\hat v}_{h}\circ {\cal
A}_{h,t}^{-1}, \;
{\hat v} \in {\cal L}^{k}(\Omega_{c}) \}, \nonumber \\
{\cal H}_{h,0}(\Omega)&  = & \{ v_{h} \, : \, \Omega
\rightarrow {\mathbb R}\,:\, v_{h}= {\hat v}_{h}\circ {\cal
A}_{h,t}^{-1}, \; {\hat v} \in {\cal L}_{0}^{k}(\Omega_{c}) \},
\end{eqnarray*}
(again, analogously to the one-dimensional setting studied in \cite[\S3.2]{macdonald:15}).
Letting ${\cal H}_{h,E_{q}}\subset {\cal L}^{k}(\Omega)$ and ${\cal
H}_{h,E_{U}}\subset {\cal L}^{k}(\Omega)$ 
be the finite dimensional approximation spaces satisfying the Dirichlet
boundary conditions for the $q_{i}$'s and $U$, respectively,
the finite element spatial discretisation of the
conservative weak formulation (\ref{eq:weak})
is therefore: find $q_{ih}(t)\in {\cal H}_{h,E_{q}}(\Omega_{t})$, 
$i=1,\ldots,5$, and 
$U_{h}\in {\cal H}_{h,E_{U}}(\Omega)$ such that $\forall v_{h} \in {\cal
H}_{h,0}(\Omega)$  
\begin{subequations}
\begin{equation}
\frac{d}{dt} \int_{\Omega} q_{ih}v_h  \; {\rm d}{{\bfa x}}-
\int_{\Omega} (\nabla \cdot ({\dot {\bfa x}}q_{ih}))\,v_h\;{\rm d}{{\bfa x}} 
= \int_{\Omega} \Gamma_{ih}\cdot \nabla v_h \; {\rm d}{{\bfa x}}- \int_{\Omega} f_{ih} v_h \; {\rm d}{{\bfa x}}, 
\quad  
\label{eq:diswcon}
\end{equation}
\begin{equation}
\int_{\Omega} {\bfa D_h}\cdot \nabla v_h \; {\rm d}{{\bfa x}} =  0.
\label{eq:dismax}
\end{equation}
\label{eq:disc}
\end{subequations}

Finally,  introducing vectors ${\bfa q}_{i}(t)$ and ${\bfa u}(t)$ which
contain the degrees of freedom defining $q_{ih}$ and $U_{h}$, respectively,
(\ref{eq:disc}) can be rewritten to obtain the highly nonlinear differential algebraic system   
\begin{subequations}
\begin{equation}
\frac{d}{dt} (M(t){\bfa q}_{i}(t)) = {\bfa G}_{i}(t,{\bfa
q}_{i}(t),{\bfa u}(t)), 
\qquad i=1,\ldots,5,
\label{eq:odes}
\end{equation}
\begin{equation}
{\bfa C}({\bfa q}_{i}(t),{\bfa u}(t))={\bfa 0}, 
\qquad i=1,\ldots,5,
\label{eq:nonsysmax}
\end{equation}
\label{eq:odesys}
\end{subequations}
where $M(t)$ is the (time-dependent) finite element mass matrix. 

\section{Derivation and discretisation of moving mesh PDEs}
\label{sec:3}
\subsection{Equations governing mesh movement}
\label{moving}
Having formulated equations to represent the physical PDEs, we now establish a mechanism 
for moving the mesh: this will be done by constructing  so-called \textit{moving mesh PDEs}.
To avoid potential mesh crossings or foldings, we derive a suitable evolution equation for the 
inverse mapping ${\cal A}_{t}^{-1}({\bfa x})={\bfa \xi}({\bfa x},t)$ rather than ${\cal A}_{t}({\bfa 
\xi})={\bfa x}({\bfa \xi},t)$ in (\ref{eq:map1}) (see, for example, the discussion in 
\cite{dvinsky:91}).  
A mesh ${\cal T}_{h,t}$ on $\Omega$ can then be generated as the pre-image of a fixed mesh ${\cal 
T}_{h,c}$ on $\Omega_{c}$.
As introduced in \cite{hr:99}, we choose the mapping ${\bfa \xi}({\bfa x})$ corresponding to a 
fixed 
value of $t$ in order to minimise the 
functional 
\begin{equation}
I[{\bfa \xi}]=\frac{1}{2} \int_{\Omega_{t}}[ (\nabla \xi)^{T} G^{-1} (\nabla \xi)+ 
(\nabla \eta)^{T} G^{-1} (\nabla \eta)]
\;{\rm d} {\bfa x},
\label{eq:var2d}
\end{equation}
where $G$ is a $2\times 2$ symmetric positive definite matrix referred to as a \textit{monitor 
matrix}, and 
$\nabla$ is the gradient operator with respect to ${\bfa x}$. Rather than directly attempt to minimise 
(\ref{eq:var2d}), a more robust procedure is to evolve the mapping according to the modified gradient flow equations 
\begin{equation}
\frac{\partial {\xi}}{\partial t}=\frac{P}{\tau} \nabla \cdot (G^{-1}\nabla{\xi}), \quad \quad {\rm and } 
\quad\quad
\frac{\partial {\eta}}{\partial t}=\frac{P}{\tau} \nabla \cdot (G^{-1}\nabla{\eta}).
\label{eq:mmpde1}
\end{equation}
Here, $\tau>0$ is a user-specified temporal smoothing parameter which affects the temporal scale 
over which the mesh adapts, and $P$ is a positive function of $({\bfa x},t)$, chosen such that 
the mesh movement has a spatially uniform time scale \cite{hr11}.

The selection of an appropriate monitor matrix is crucial to the success of mesh adaptation. 
In this paper, we will consider the monitor matrix proposed by Winslow \cite{winslow:67} 
\begin{equation}
G=\left [ \begin{array}{cc}
w& 0 \\
0 & w \end{array} \right ],
\label{eq:winmon}
\end{equation}
where $w({\bfa x},t)$ is a positive weight function called a {\em monitor function}. The choice of 
the monitor function should ideally be based on a local a posteriori error estimate but if no such 
estimate exists then the monitor function can be any smooth function designed to adapt the mesh 
towards important solution features. Suitable choices for $w$ for applications to the ${\bfa 
Q}$-tensor equations are discussed below.

In practice, we interchange the roles of the dependent 
and independent variables in (\ref{eq:mmpde1}),  
since it's the location of the physical mesh points $\{{\bfa x}_{i}(t)\}_{i=1}^{\cal N}$ that defines the mapping ${\cal A}_{t}$. 
With a Winslow-type monitor matrix (\ref{eq:winmon}) the resulting MMPDEs take the form  
\begin{equation}
\tau \frac{\partial {\bfa x}}{\partial t}=
P( a{\bfa x}_{\xi\xi}+b{\bfa x}_{\xi\eta}+c{\bfa x}_{\eta\eta}+d{\bfa x}_{\xi}+
e{\bfa x}_{\eta}), \quad (\xi,\eta)\in \Omega_{c},
\label{eq:mmpde}
\end{equation}
where
\[
a=\frac{1}{w} \frac{x_{\eta}^{2}+y_{\eta}^{2}}{J^{2}}, \quad 
b=-\frac{2}{w} \frac{(x_{\xi}x_{\eta}+y_{\xi}y_{\eta})}{J^{2}},\quad
c=\frac{1}{w} \frac{x_{\xi}^{2}+y_{\xi}^{2}}{J^{2}},
\]
\[
d=\frac{1}{(w J)^{2}}[ w_{\xi}(x_{\eta}^{2}+y_{\eta}^{2})-
                                            w_{\eta}(x_{\xi}x_{\eta}+y_{\xi}y_{\eta}),
\]                                           
\[                                            
e=\frac{1}{(w J)^{2}}\left [
 -w_{\xi}(x_{\xi}x_{\eta}+y_{\xi}y_{\eta})+w_{\eta}(x_{\xi}^{2}+y_{\xi}^{2})\right ],
 \] 
 and $J=x_{\xi}y_{\eta}-x_{\eta}y_{\xi}$ is the Jacobian of ${\cal A}_{t}$. To complete 
the specification of the coordinate transformation, the MMPDE (\ref{eq:mmpde}) must 
be supplemented by suitable boundary conditions 
${\bfa g}({\bfa \xi},t)$, ${\bfa \xi}\in \partial \Omega_{c}$; these are obtained using 
a one-dimensional moving mesh approach.

The numerical solution of (\ref{eq:mmpde}) requires spatial and temporal discretisation. 
We discretise in space using standard linear Galerkin finite elements. For time discretisation, we
use a backward Euler integration scheme to update the solution at $t=t^{n+1}$ and, to avoid 
solving nonlinear algebraic systems, we evaluate the coefficients $a,c,\ldots,e$ at the time $t=t^{n}$. 
We therefore seek 
${\bfa x}_{h}^{n+1}\in ({\cal L}^{1}(\Omega_{c}))^{2}$ 
such that 
\begin{eqnarray}
\tau\int_{\Omega_{c}} \left( \frac{{\bfa x}_{h}^{n+1}-{\bfa x}_{h}^{n} }{\Delta t}\right ) \cdot {\bfa {\hat v}}_{h}\; {\rm d}{\bfa \xi}+ 
\int_{\Omega_{c}} \left [ ({\bfa x}_{h}^{n+1})_{\xi}\cdot (a^{n}{\bfa {\hat v}}_{h})_{\xi} + ({\bfa x}_{h}^{n+1})_{\eta}\cdot(c^{n}{\bfa {\hat v}}_{h})_{\eta} \right .& & \nonumber \\
& &  \hspace{-8.5cm}+\tfrac{1}{2}[({\bfa x}_{h}^{n+1})_{\xi}\cdot(b^{n}{\bfa {\hat v}}_{h})_{\eta}+({\bfa x}_{h}^{n+1})_{\eta}\cdot (b^{n}{\bfa {\hat v}}_{h})_{\xi}]\nonumber  \\
& & \left . \hspace{-7cm} 
-[d^{n}({\bfa x}_{h}^{n+1})_{\xi}+e^{n}({\bfa x}_{h}^{n+1})_{\eta}]\cdot {\bfa {\hat v}}_{h} \right ] \:{\rm d} {\bfa \xi} = 0,
\label{eq:wkmmpde}
\end{eqnarray}
for all ${\bfa {\hat v}}_{h}\in ({\cal L}^1_{0}(\Omega_{c}))^{2}$. The resulting linear systems are 
solved using the iterative method BiCGSTAB \cite{vorst92} with an incomplete LU (ILU) 
factorization \cite{vorst77} as a preconditioner. An analysis of the performance of this 
iterative solver for the discretised MMPDE equations can be found in \cite{beckett:02}.

\subsection{Details of the monitor functions}
\label{subsectmon}
An appropriate choice of monitor function $w({\bfa x},t)$ in (\ref{eq:winmon}) is essential to the 
success of any adaptive moving mesh 
method. In this paper we consider two-dimensional analogues of the monitor functions which were 
shown in \cite{macdonald:15} to be appropriate for one-dimensional $\bfa Q$-tensor models.
We first describe these assuming that we have an input function $\mathcal{T}({\bfa x},t)$ which 
represents a physical quantity associated with the particular problem under 
consideration: a discussion of appropriate input functions for our problem involving 
finite element approximation of the ${\bfa Q}$-tensor matrix 
follows in \S\ref{input}.

We consider three different forms of monitor function: 
\begin{itemize}
 \item \textbf{AL}. This is based on a measure of the arc-length of $\mathcal{T}$:
\begin{equation}
w(\mathcal{T}({\bfa x},t))= \left({1+\left|\nabla \mathcal{T}({\bfa 
x},t)\right|^2}\right)^{\frac{1}{2}}.
\label{arc}
\end{equation}
 \item \textbf{BM1}. This is a generalisation of the Beckett-Mackenzie monitor function 
introduce in \cite{bm00a,bm01b}, based on first-order partial derivatives of $\mathcal{T}$:
\begin{equation}
w(\mathcal{T}({\bfa x},t))= \alpha({\bfa x},t) + \left|\nabla \mathcal{T}({\bfa x},t), 
\right|^{\frac{1}{m}}
\label{bm1}
\end{equation}
where scaling parameters $\alpha$ and $m$ are discussed below. 
 \item \textbf{BM2}.
 This is a second variant of the Beckett-Mackenzie monitor function which takes into account 
information about second-order partial derivatives of $\mathcal{T}$: 
\begin{equation}
w(\mathcal{T}({\bfa x},t))= \alpha({\bfa x},t) + 
\left( \sqrt{\left(\frac{\partial^2 \mathcal{T}}{\partial x^2}\right)^2
+2\left(\frac{\partial^2 \mathcal{T}}{\partial x \partial y}\right)^2 + \left(\frac{\partial^2 
\mathcal{T}}{\partial y^2}\right)^2}
\right)^{\frac{1}{m}}
\label{bm2}
\end{equation}
again involving scaling parameters $\alpha$ and $m$.
\end{itemize}

The value of the parameter $\alpha$ in (\ref{bm1}) and (\ref{bm2}) 
is determined \textit{a posteriori} from the numerical approximation itself as  
\begin{displaymath}
\alpha({\bfa x},t)=\textrm{max}\left\{1,\;\frac{1}{\textrm{meas}(\Omega)}\int_{\Omega} 
\mathcal{I}\;
\textrm{d}{\bfa x}\right\},
\label{alpha}
\end{displaymath}
where
\begin{displaymath}
\mathcal{I}=\left|\nabla \mathcal{T}({\bfa x},t) \right|^{\frac{1}{m}} 
\end{displaymath}
for BM1 and 
\begin{displaymath}
\mathcal{I}=
\left( \sqrt{\left(\frac{\partial^2 \mathcal{T}}{\partial x^2}\right)^2
+2\left(\frac{\partial^2 \mathcal{T}}{\partial x \partial y}\right)^2 + \left(\frac{\partial^2 
\mathcal{T}}{\partial y^2}\right)^2}
\right)^{\frac{1}{m}}
\end{displaymath}
for BM2. Its purpose is to avoid mesh starvation external to layers in the solution as, without it, 
the resulting mesh would have almost all mesh points clustered inside the layers due to the monitor 
function being very small elsewhere. The lower bound on $\alpha$ in (\ref{alpha})
removes unwanted oscillations in the mesh trajectories caused by amplification of
errors in approximating $\mathcal{T}(\bfa{x},t)$ which could otherwise cause 
the mesh to adapt incorrectly.  
The parameter $m$ in (\ref{bm1}) and (\ref{bm2}) also helps to regulate mesh clustering:
when $m>1$, any large variations in $\mathcal{T}({\bfa{x}},t)$ are smoothed  so that mesh 
points are more evenly distributed over the domain.  
In \cite{bm01b} it is shown that, for a function in one dimension with a boundary layer, the
optimal rate of approximation order using polynomial elements 
of degree $p$ can be obtained by ensuring that the parameter $m\geq p+1$.  With no
specific guidance on the choice of $m$ in higher-dimensional settings, we follow our work
in \cite{macdonald:15} and choose $m=3$.

In general, the monitor function often has large spatial and temporal variations, so to improve the 
robustness of the moving mesh method we employ both spatial and temporal smoothing procedures. 
This results in an MMPDE that is easier to integrate forward in time and a smoother mesh, which can 
improve spatial solution accuracy.  Temporal smoothing is done by under-relaxing the monitor 
function so that the monitor function at the current time level $n$ is given by
\begin{equation}
w^{n}=(1-\omega)w^{n}+\omega w^{n-1},
\end{equation}
where $0<\omega<1$ is an under-relaxation parameter (in this paper, we set $\omega=0.8$).
Following \cite{hr11}, spatial smoothing of the monitor function is done by taking a local 
average
of the monitor function across cells that have a common vertex. That is, the smoothed monitor 
function $\tilde{w}$ is defined as
\begin{equation}
\tilde{w}({\bfa x}_m)=\frac{\int_{{C{({\bfa \xi}_m)}}}w({\bfa x}({\bfa \xi}))\textrm{d}{\bfa 
\xi}}{\int_{{C{({\bfa \xi}_m)}}}\textrm{d}{\bfa \xi}},
\end{equation}
where ${\bfa x}_m\in\Omega$ is a mesh point in the physical domain, ${\bfa \xi}_m\in\Omega_{c}$ is 
the corresponding a mesh point in
the computational domain, and $C{({\bfa \xi}_m)}\subset \Omega_{c}$ represents all neighbouring 
cells of vertex
${\bfa \xi}_m$.  If required, spatial smoothing can be repeated in an iterative fashion to
further smooth the monitor function.

\subsection{Details of the monitor input functions}
\label{input}
Having specified monitor functions, it remains to decide on an 
appropriate input function $\mathcal{T}(\bfa{x},t)$. 
We will consider two variants:
\begin{enumerate}
\item[\textbf{a.}] \textbf{Order parameter}. We recall from \S\ref{sec:qtensor} that, for a 
uniaxial
state with scalar order parameter $S$, we have $S^2=\mathrm{tr}({\bfa Q}^2)$.
This has led to the function
\begin{equation}
\mathcal{T}(\bfa{x},t)=\mathrm{tr}(\bfa{Q}^2),
\label{trace}
\end{equation}
being used to generate monitor functions in previous studies 
\cite{abl10,abl11,rn07,rn08}. This quantity is known to vary rapidly in regions where order 
reconstruction occurs, and was shown in \cite{macdonald:11} to be ideal for certain 
one-dimensional uniaxial problems.

\item[\textbf{b.}] \textbf{Biaxiality}. For biaxial problems, an alternative input function 
(based on the direct invariant measure of biaxiality used in \cite{barberi:04}) is
\begin{equation}
\mathcal{T}(\bfa{x},t)=
\left[1-\frac{6\,{\textrm{tr}}(\bfa{Q}^3)^2}{{\textrm{tr}}(\bfa{Q}
^2)^3 } \right]
^{\frac{1}{2}}.
\label{biax_b}
\end{equation}
This takes values ranging from 0 (for a uniaxial state) to 1 (for a wholly biaxial state).
\end{enumerate}
In the numerical experiments in \S\ref{results}, we compare the performance of the AL, BM1 and BM2 
monitor functions with various input functions. Details of the specific combinations highlighted in 
the results presented are summarised in Table~\ref{montab}.
\begin{table}[htb]
\begin{center}
\begin{tabular}{|c|c|c|c|c|}
\hline
Method&\textbf{AL}&\textbf{BM1a}&\textbf{BM1b}&\textbf{BM2b}\\
\hline
Monitor function& (\ref{arc})& (\ref{bm1})& (\ref{bm1})& (\ref{bm2})\\
\hline
Input function& (\ref{trace})& (\ref{trace})& (\ref{biax_b})& (\ref{biax_b})\\
\hline
\end{tabular}
\end{center}
\caption{Details of monitor function construction. \label{montab}}
\end{table}

\section{Solution algorithm}
\label{algorithm}
The numerical algorithm for solving the full problem involves an iterative solution 
strategy as originally proposed in \cite{beckett:02,beckett:01}. Full details of how 
it can be used in a ${\bfa Q}$-tensor setting are given in \cite[\S5]{macdonald:15}, so are not 
reproduced here. Instead, we simply highlight the main features and point out any modifications 
needed for application in two dimensions.

The iterative solution algorithm involves completely decoupling the solution of the physical PDE 
system (\ref{eq:odesys}) from the solution of the MMPDE (\ref{eq:mmpde}). This has a key advantage 
in that different convergence criteria can be used for the two systems: it is well understood that 
the computational mesh rarely needs to be resolved to the same degree of accuracy as the 
solution of the physical PDEs. The system for the ${\bfa Q}$-tensor components (\ref{eq:odes}) is 
integrated forward in time using a second-order singly diagonally implicit Runge-Kutta (SDIRK2) method, 
with the electric potential values being updated by solving (\ref{eq:nonsysmax}) at each step of 
the Newton iteration used to generate intermediate stages (see \cite[\S5.2]{macdonald:15}).
We also use an adaptive time-stepping procedure when integrating forward in time. This is 
based on both the computed solutions of (\ref{eq:odes}) for ${\bfa q}_i$ and on the solution of the 
MMPDE (\ref{eq:mmpde}). The specific two-dimensional error indicator used here is
\begin{eqnarray*}
\textrm{E}_i=\left(\sum_{j=1}^{N_{\textrm{E}}}\textrm{area}_{\Delta^n_j}\left({e^{n+1}_{i,j}}
\right)^2\right)^{\frac{1}{2}},
\end{eqnarray*}
where $N_{\textrm{E}}$ is the number of mesh elements and $\Delta^n_{j}$ denotes the 
$j^{\textrm{th}}$ element of the mesh at time level $n$. The individual error terms
\begin{eqnarray*}
e^{n+1}_{i,j}={q}^{n+1}_{i,j}-{{\hat q}}^{n+1}_{i,j}
\end{eqnarray*}
represent the errors at the midpoint of each triangular element (cf.\ 
\cite[\S5.3]{macdonald:15}).  They are calculated at time level $n+1$ using the $j^{\textrm{th}}$ 
entries of the solution vector
${q}^{n+1}_{i,j}$ and the vector ${\hat q}^{n+1}_{i,j}$, which is the embedded first-order SDIRK 
approximation to ${q}^{n+1}_{i,j}$. Note that this last vector can be 
obtained from the SDIRK2 scheme at no extra computational cost. For full details of how this 
indicator is subsequently used to develop an adaptive time-stepping strategy, see 
\cite[\S5.3]{macdonald:15}.

\section{Numerical results}
\label{results}
In this section, we illustrate the performance of our novel time-adaptive 
method with a moving two-dimensional finite element mesh using two test problems involving ${\bfa 
Q}$-tensor models of liquid crystal cells. Note that we use quadratic finite elements on triangualr 
meshes throughout.

\subsection{Test problem 1: resolving the core structure of a stationary defect}
\label{static}
\begin{figure}[!htb]
\begin{center}
 \leavevmode
 \mbox{
\begin{minipage}{2.3in}
\hspace*{-2cm}
\scalebox{0.4}{\includegraphics{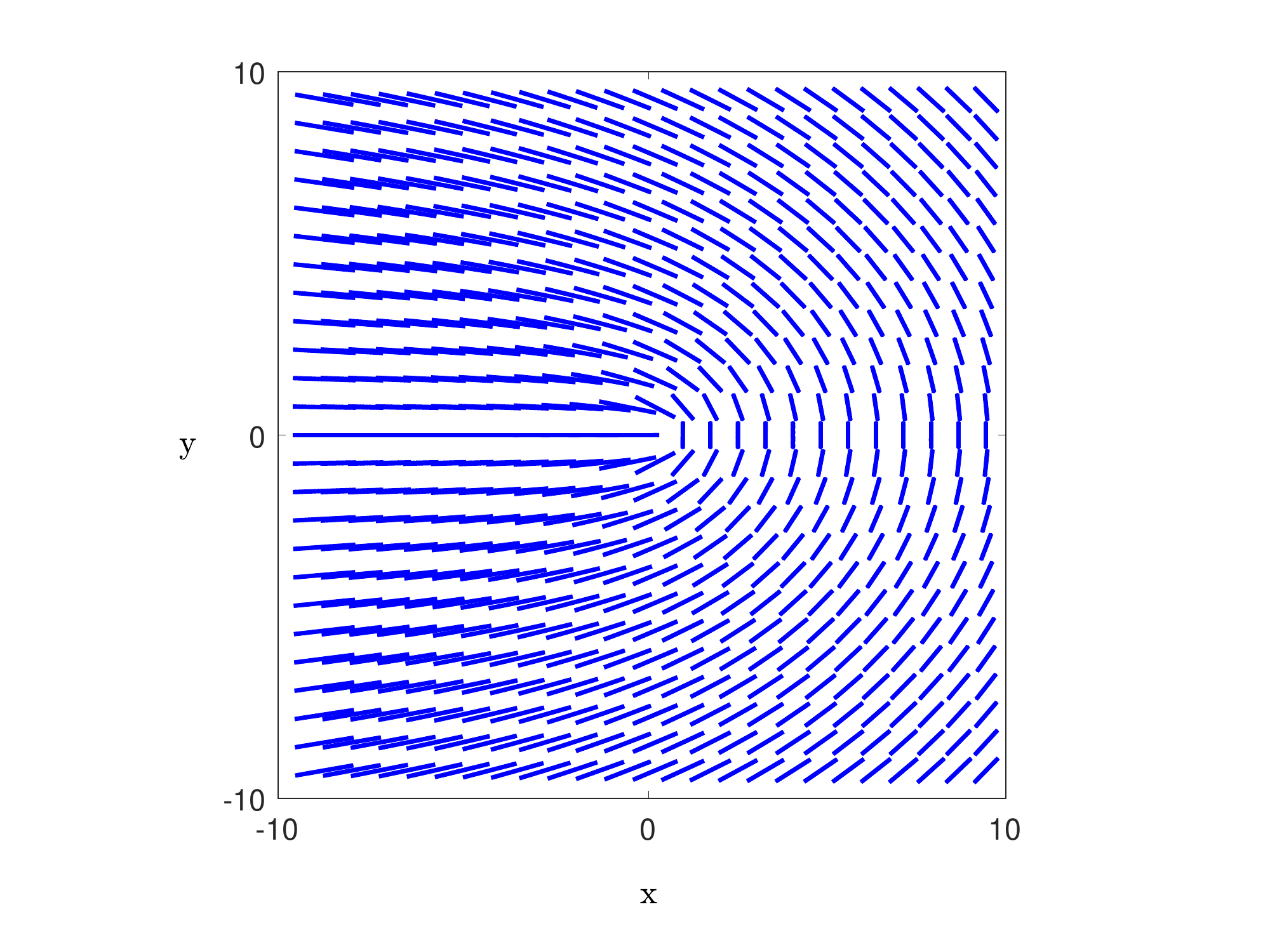}}\\
 {\small (a) Director field.}\end{minipage}
\qquad
 \begin{minipage}{2.3in}
\hspace*{-2cm}
\scalebox{0.37}{\includegraphics{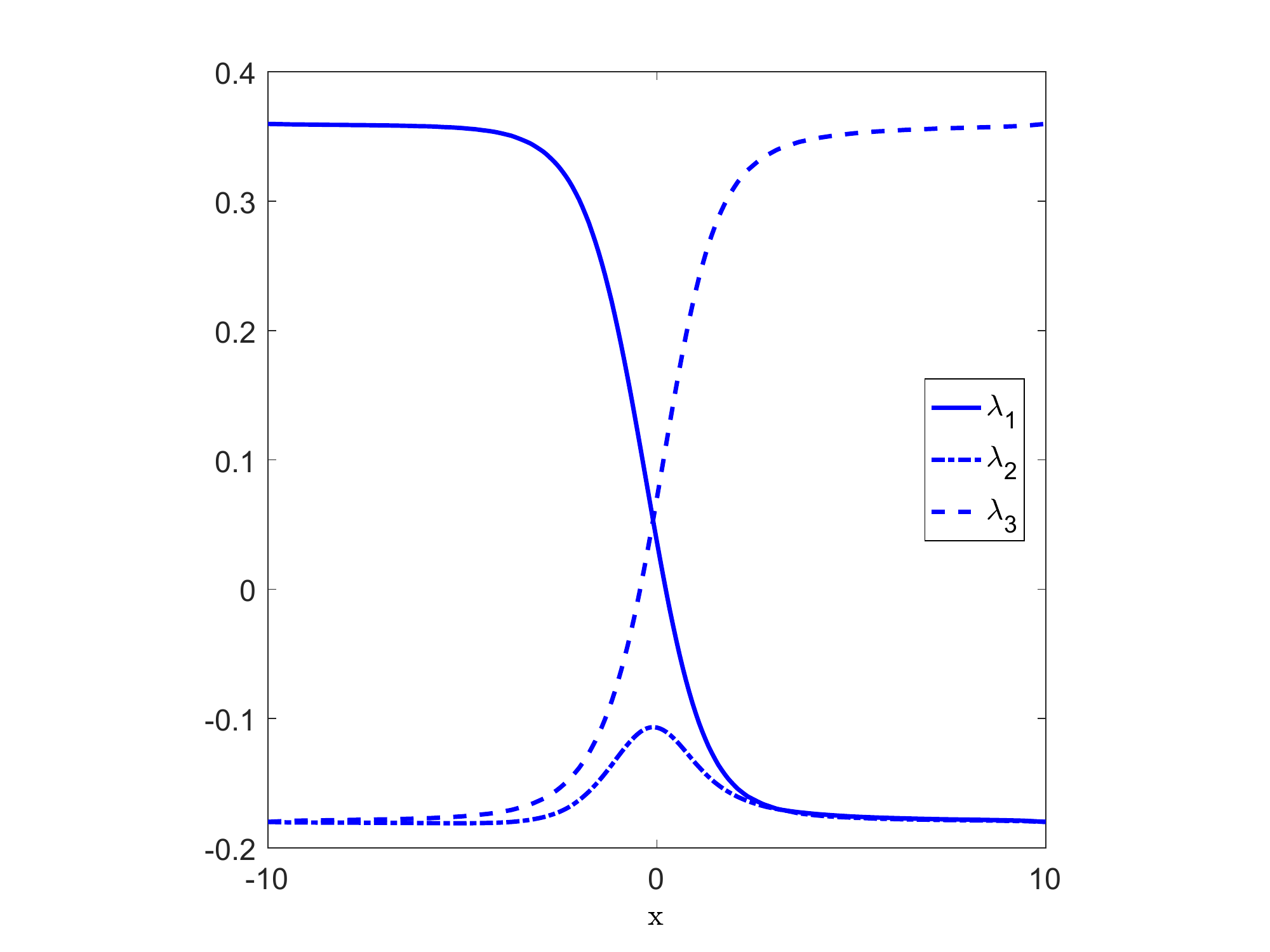}}\\
 {\small (b) Eigenvalues of ${\bfa Q}$.
}\end{minipage} 
}
\end{center}
\caption{Plots of the director field and eigenvalues exchange along the line 
$y=0$ for a +1/2 defect.
\label{fig:dirseigs}}
\end{figure}

We begin by considering the resolution of a stationary liquid crystal defect, as this problem is 
ideal for assessing the ability of particular monitor functions to resolve the core defect 
structure. Specifically, we assume we have a disclination line in the $z$-direction and 
that far from the defect core the director $\bfa{n}$ lies in the $x$-$y$ plane.  
Such a defect can be simulated by imposing the initial condition on the director 
\begin{equation}
{\bfa n} = \left(\cos(d_{i}\theta),\sin(d_{i}\theta),0\right), \quad 0\leq \theta \leq 2\pi, 
\label{eq:dirdefect}
\end{equation}
where $d_{i}$ is referred to as the {\em disclination index}. We note that on travelling a closed 
path around the 
disclination line, the director rotates through the angle $2\pi d_{i}$. Here we restrict 
our attention to the 
case $d_{i}=1/2$:
a plot of the steady state director field is shown in Figure \ref{fig:dirseigs}(a), where we 
observe that the director rotates through an angle of $\pi$ radians as the centre of 
the defect is circled. 
Using numerical simulations of a ${\bfa Q}$-tensor model, Schopohl \& 
Sluckin \cite{schopohl:87} found that the defect core structure was contained in a circular 
region of radius approximately $5\zeta$, where $\zeta$ is the nematic coherence length. For this 
example, we solve the ${\bfa Q}$-tensor equations on the square region 
$[-10\zeta,10\zeta]^{2}$, with
$\zeta\approx4.06$ nm (as commensurate with our use of the physical parameters in 
\cite{barberi:04}).
 On the domain boundaries, we impose Dirichlet conditions corresponding to 
the director being given by (\ref{eq:dirdefect}) and the order parameter $S$ taking a value 
associated with the equilibrium (nematic) phase ($S=S_{eq}\simeq 0.65$ with our parameters). An 
analytical 
solution of the ${\bfa Q}$-tensor equations does not exist for this problem, so a reference 
solution was obtained using a 
fine adaptive mesh using the BM2b monitor function with $N=5334$ quadratic triangular elements,
using a time step $\Delta t=10^{-8}$ until a steady state solution was 
obtained. We consider this reference solution to be independent of the specific 
choice of the monitor function as calculations using reference solutions based on the other three 
monitor functions gave very similar results. 

\begin{figure}[!htb]
\begin{center}
 \leavevmode
 \mbox{
\begin{minipage}{2.3in}
\hspace*{-1cm}
\scalebox{0.33}{\includegraphics{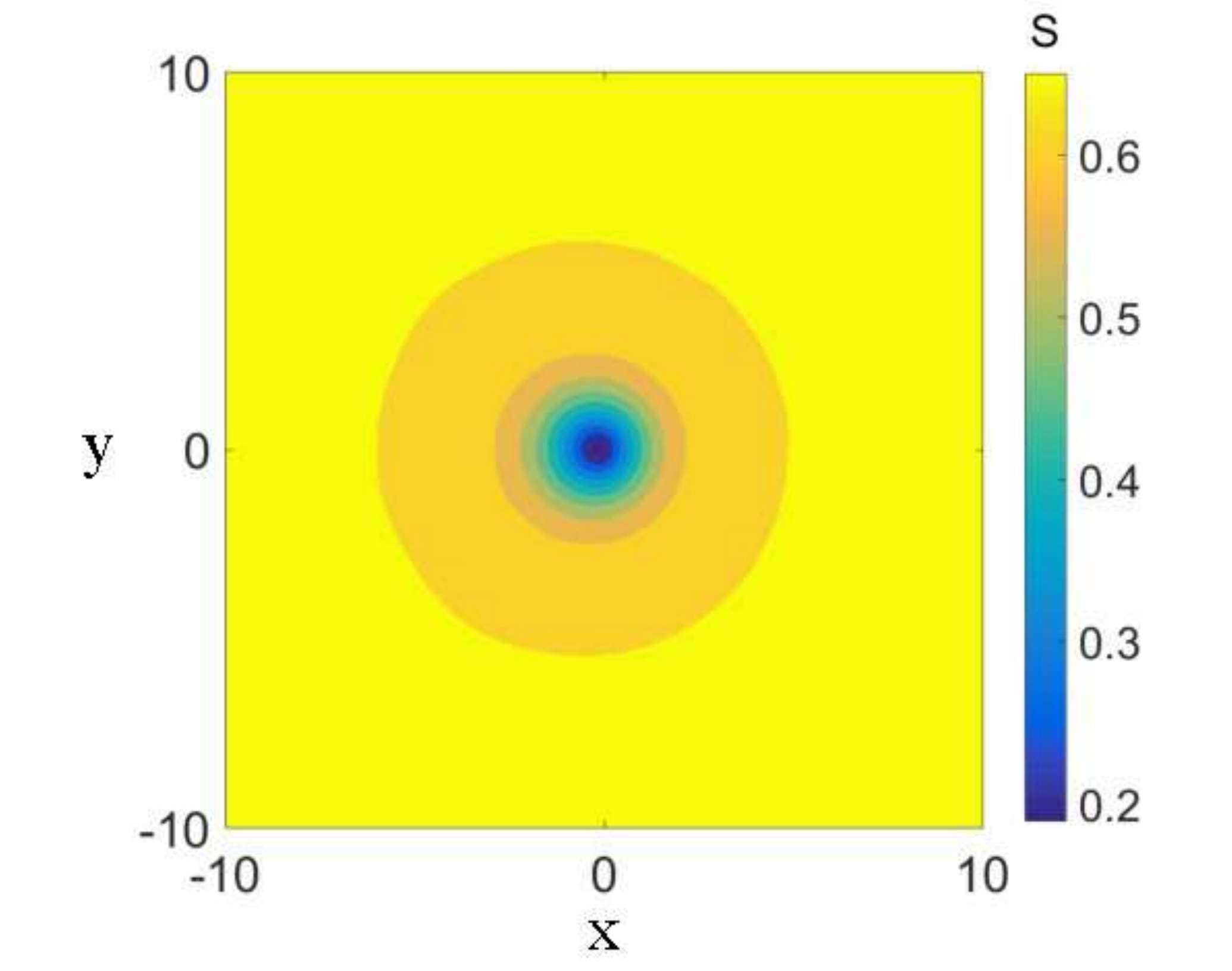}}\\
 {\small (a) Order parameter.}\end{minipage}
 \quad
 \quad
 \begin{minipage}{2.3in}
\hspace*{-1cm}
\scalebox{0.32}{\includegraphics{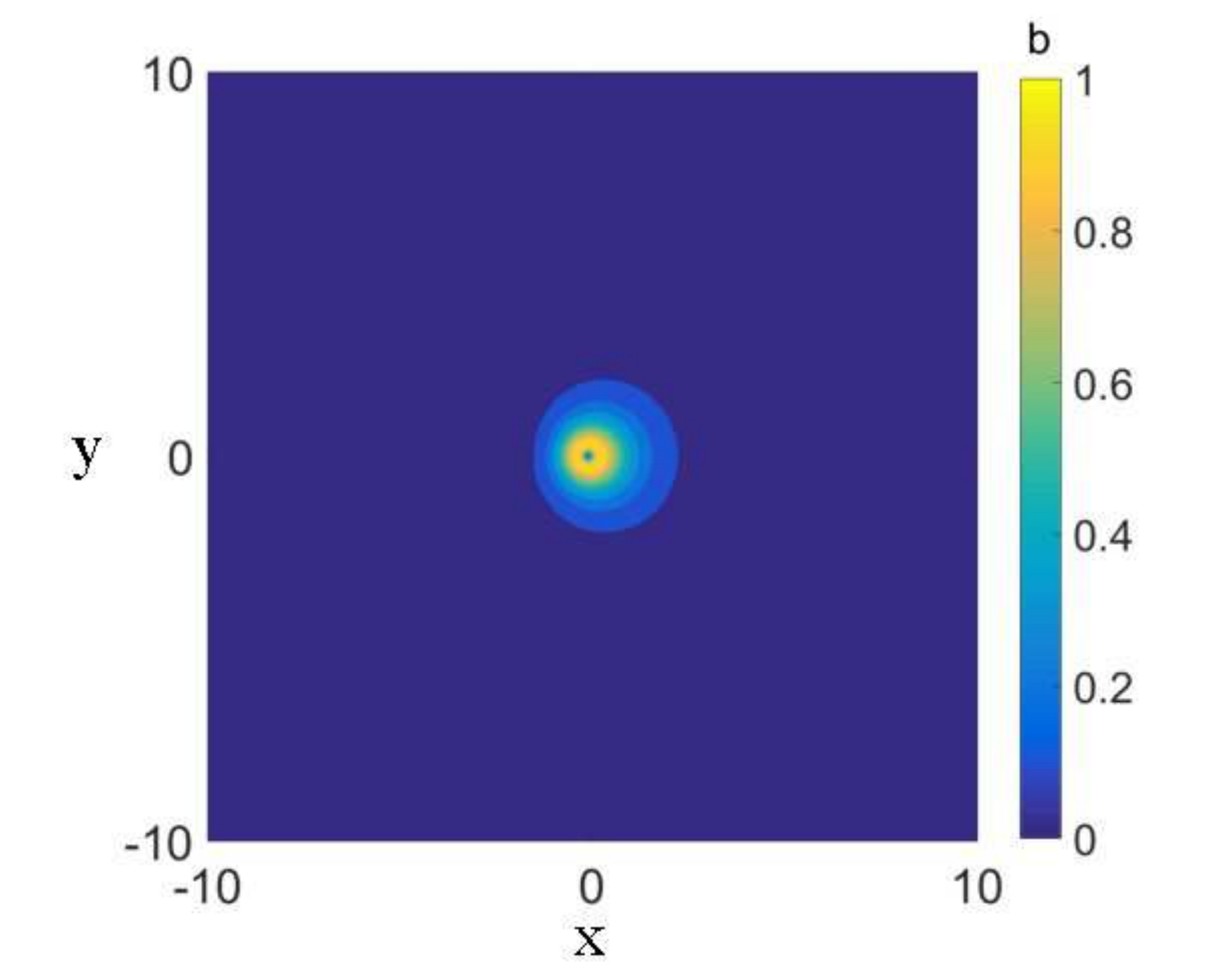}}
 {\small (b) Biaxiality.
}\end{minipage} 
}
\end{center}
\caption{Contour plots of the order and biaxiality for +1/2 defect.
 \label{Sb_d1}}
\end{figure}

Figure \ref{fig:dirseigs}(b) shows the three (numerically computed) eigenvalues of the 
${\bfa Q}$-tensor along the line $y=0$ at the centre of the cell. 
We observe that an exchange of eigenvalues takes place at the centre of the core region, as the 
material passes through a biaxial transition (corresponding 
to the switch from horizontally to vertically aligned directors along $y=0$ 
in Figure~\ref{fig:dirseigs}(a).
Contour plots of the computed order parameter $S=tr({\bfa Q}^{2})^{1/2}$ (cf.\ (\ref{trace})) and 
the 
biaxiality (as measured by (\ref{biax_b})) are shown in Figure \ref{Sb_d1}.
We can see that the order parameter takes its equilibrium value $S\approx 0.65$ outside 
a central circular region of diameter $10\zeta$.  Furthermore, within this region the order 
parameter varies significantly within a core of diameter $2\zeta$. We can also see that, outside a 
core of diameter of approximately $4\zeta$, the biaxiality is zero, and inside it has a 
volcano-like structure with a rim of value 1 representing the purely biaxial state, and a base 
with value 0 representing the uniaxial state: this sudden variation takes place over a core only a 
few nanometres in diameter, so is very difficult to capture accurately using a numerical method.  

\subsubsection{Estimated rate of spatial convergence} 
\label{conv}
It is important to check that approximations obtained on a sequence 
of adaptive meshes are convergent as we increase the number of mesh 
elements. In \cite{macdonald:11}, we presented convergence results for a 
scalar model of a one-dimensional uniaxial problem, with a much more 
complicated one-dimensional Pi-cell order reconstruction problem, modelled 
by a full ${\bfa Q}$-tensor model, being considered in \cite{macdonald:15}.  
In a similar vein, we now consider the rate of spatial convergence of the 
moving mesh finite element approximation for this fully two-dimensional defect 
model problem. 

In what follows, we use $q_{i\ast}({\bfa x},t)$ to denote the reference 
approximation to the exact solution $q_i({\bfa x},t)$, and $q_{iN}$ to denote the finite element 
approximation calculated on a grid with $N$ quadratic elements.
We assume throughout that
\bd
|q_{i\ast}({\bfa x},t)-q_i({\bfa x},t)|\ll \left|q_{i\ast}({\bfa x},t)-q_{iN}({\bfa x},t)\right|.
\ed
The error in the approximation $q_{iN}$ is denoted by $e^N_{q_i}$. 
Since the approximate solution grid points will not in general coincide with the
reference grid points, the reference solution at a fixed time $t=t^\ast$, $q_{i\ast}({\bfa 
x}_{jk},t^\ast)$, is interpolated (using the MATLAB function \texttt{TriScatteredInterp} 
\cite{matlab}) onto the approximate solution grid. The spatial error in the $l_\infty$ norm is 
then estimated using the maximum error computed at the grid nodes, that is,
\be
\|{\bfa
e}_{q_i}^N\|_{l_\infty}=\max_{j,k=0,\ldots,N}|q_{i\ast}({\bfa x}_{jk},t^\ast)-q_{iN}({\bfa x}_{jk},
t^\ast)|.
\label{linf}
\ee

\begin{figure}[!htb]
\begin{center}
 \leavevmode
 \mbox{
\begin{minipage}{2.3in}
\hspace*{-1cm}
\scalebox{.35}{\includegraphics{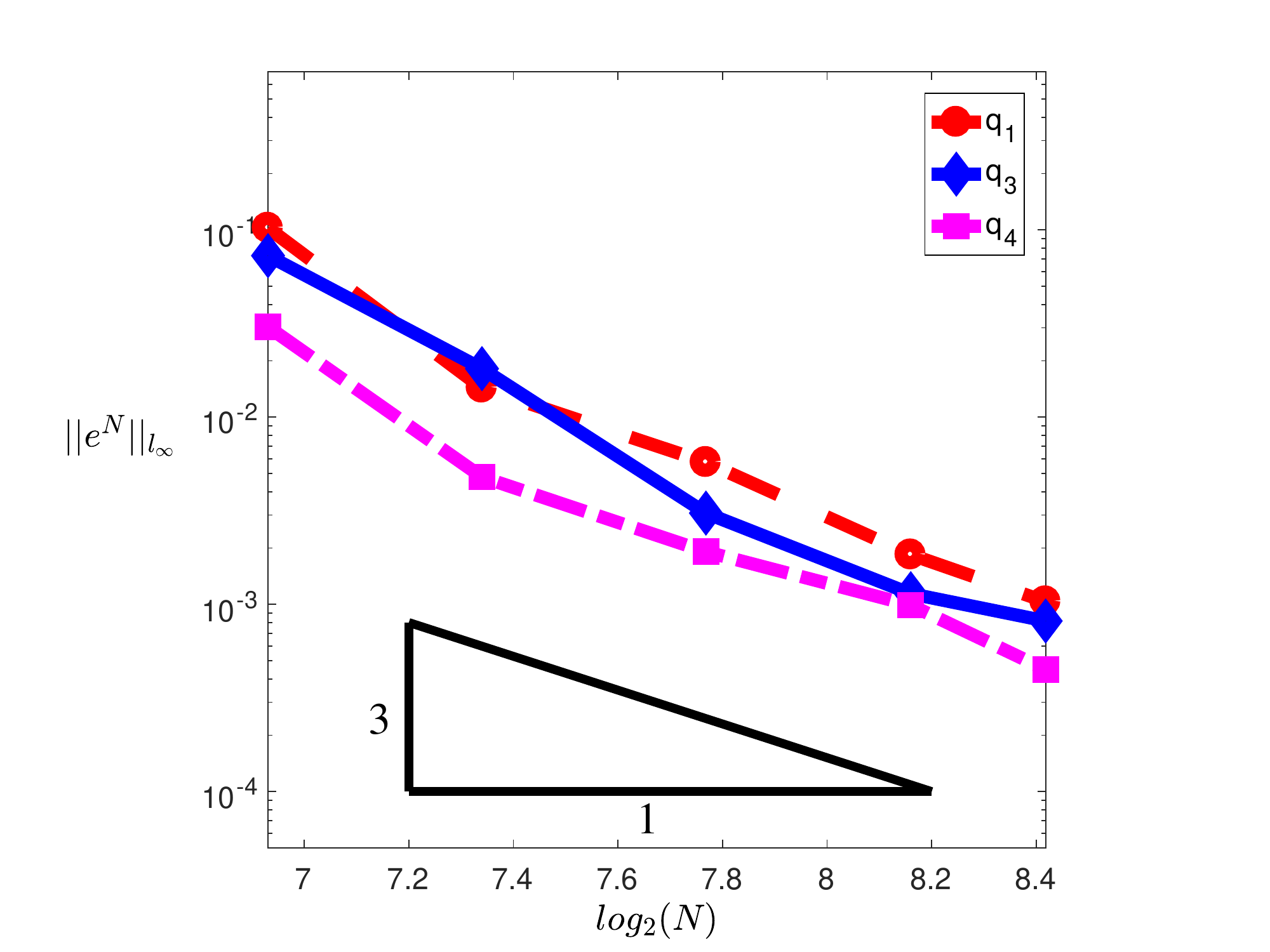}}\\
 {\small (a) $l_\infty$ error (\ref{linf}) for $q_1,q_3,q_4$.}\end{minipage}
 \begin{minipage}{2.6in}
\hspace*{-1cm}
\scalebox{0.3}{\includegraphics{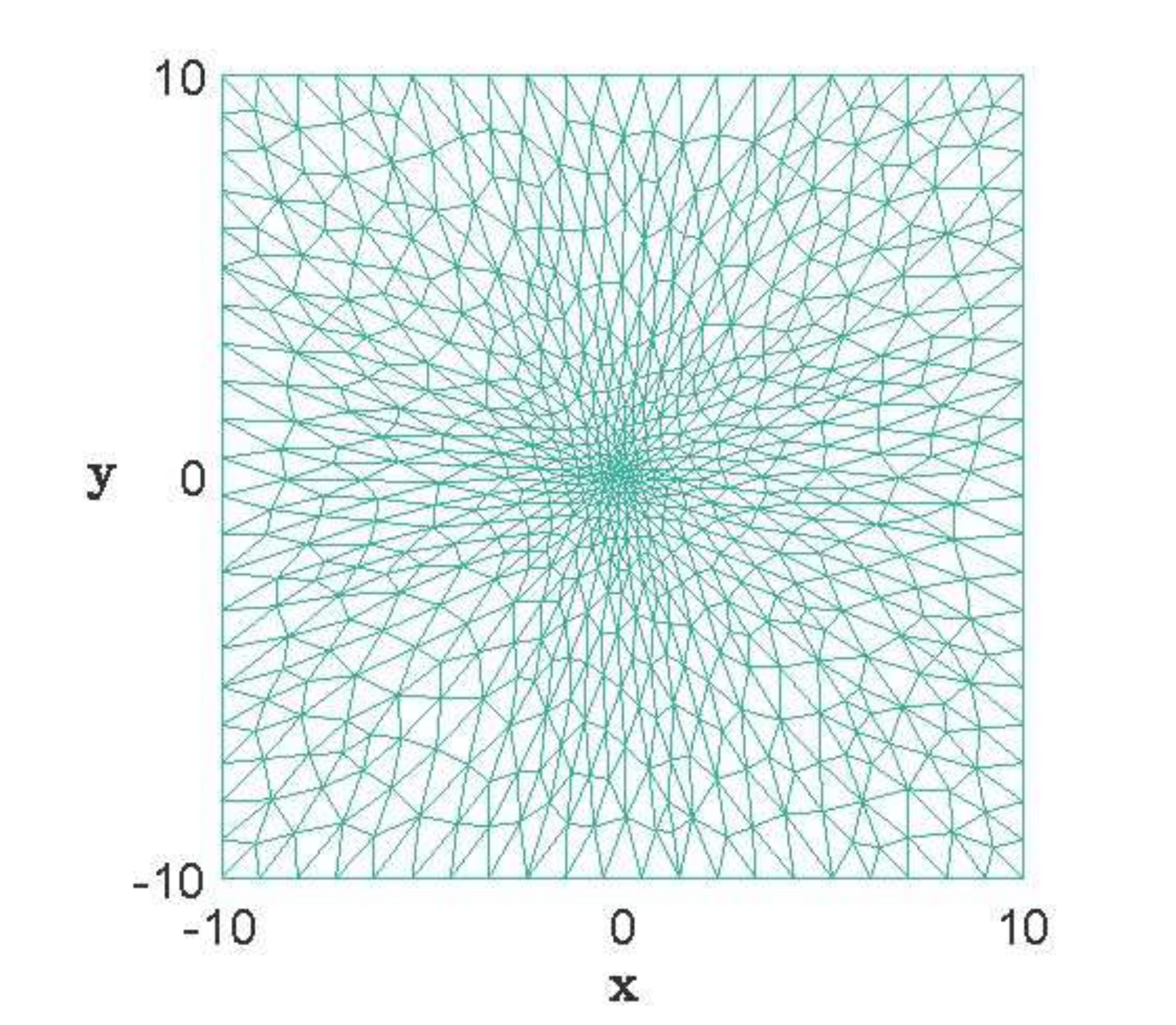}}\\
 {\small (b) Adapted mesh with $N=1388$ quadratic elements.
}\end{minipage} 
}
\end{center}
\caption{Error norm and sample adapted grid for +1/2 defect. \label{err_grid}}
\end{figure}

We fix time $t^\ast=0.2$ ms as by this time the solution has entered a steady state.
All approximate solutions are obtained using the BM2b monitor function. 
The error norm (\ref{linf}) for the non-zero components
of $\bfa{Q}$ (components $q_2$ and $q_5$ are exactly zero for this problem) for the +1/2 defect 
problem are plotted in Figure \ref{err_grid}(a) for various values of $N$. 
We observe that $\|e^N_{q_i}\|_{l_\infty}$ appears to converge at the rate $\mathcal{O}(N^{-3})$ 
which is the optimal rate expected using quadratic triangular elements.

\subsubsection{Resolving the defect core}

\begin{figure}[ht]
\scalebox{.34}{\includegraphics{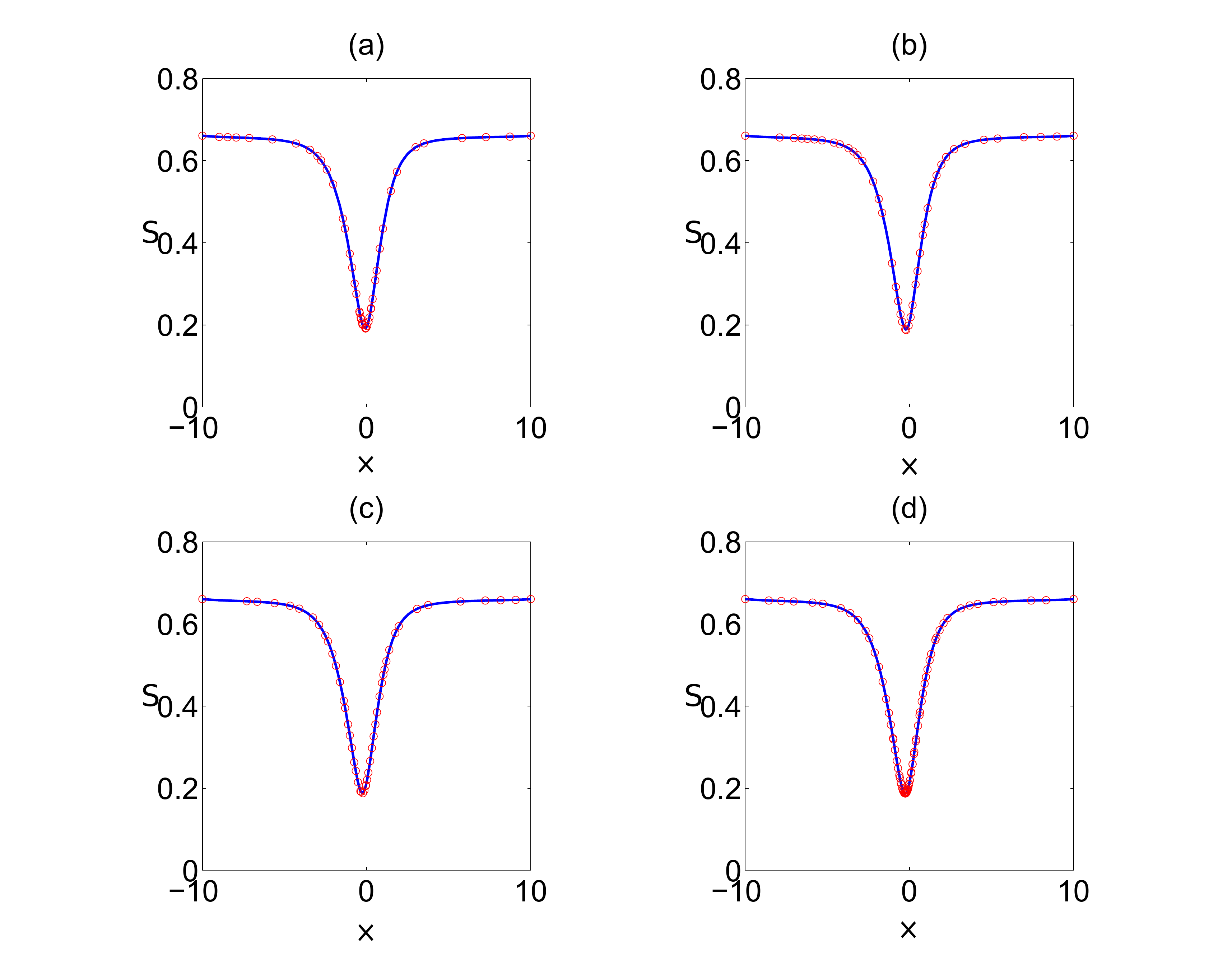}}
\caption{Order parameter along the cross section $y=0$ for the +1/2 defect, obtained 
using 1946 quadratic elements using the monitor functions: (a) AL; (b) BM1a; (c) BM1b; (d) 
BM2b. \label{opplot}}
\end{figure}

\begin{figure}[!ht]
\begin{center}
\scalebox{.6}{\includegraphics{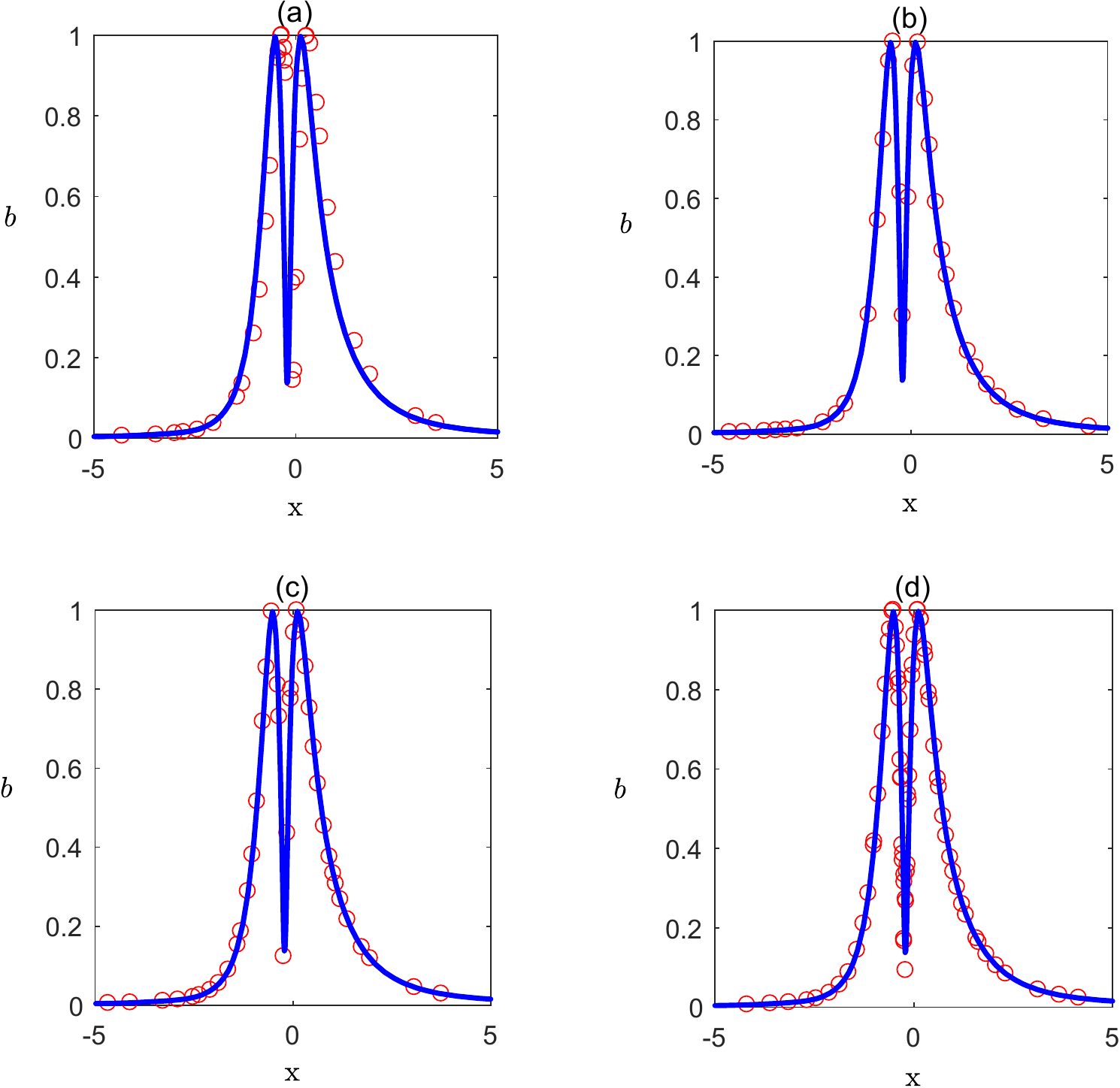}}
\caption{Biaxiality along the cross section $y=0$ for the +1/2 defect, obtained 
using 1946
quadratic elements using the monitor functions: (a) AL; (b) BM1a; (c) BM1b; (d) 
BM2b. \label{biax}}
\end{center}
\end{figure}

An example of an adapted mesh using $N=1388$ elements and the BM2b monitor function is shown in 
Figure~\ref{err_grid}(b). Although it is clear that the mesh has been adapted isotropically towards 
the core of the defect, at this scale it is difficult to observe any detail of exactly how the
adaptivity resolves the defect core structure. To give some insight into the resolution 
obtained using the different monitor functions, Figures~\ref{opplot} and \ref{biax} show
cross-sections (taken along the line $y=0$)  of the order parameter and biaxiality, respectively 
(plotted as solid blue lines). 
The location of grid nodes for the monitor and input function combinations listed in 
Table~\ref{montab} are also plotted (as red circles) in each case to help visualise
how each method copes with adapting to the nano-structure of the defect core. 
All of the meshes are clearly adapting to resolve the core structure of
both the order parameter and the biaxiality. However, we note that the BM2b monitor function in 
particular has placed a significant number of nodes exactly at the defect core, 
right inside the volcano structure coming from the biaxiality. We know from our previous
experience with the one-dimensional Pi-cell problem \cite{macdonald:15} that it is particularly 
difficult to resolve this structure, but the BM2b monitor function is still doing a good job here 
for the full two-dimensional problem.

\begin{figure}[!htb]
 \leavevmode
 \mbox{
 \hspace*{-2cm}
\begin{minipage}{18cm}
\scalebox{0.2}{\includegraphics{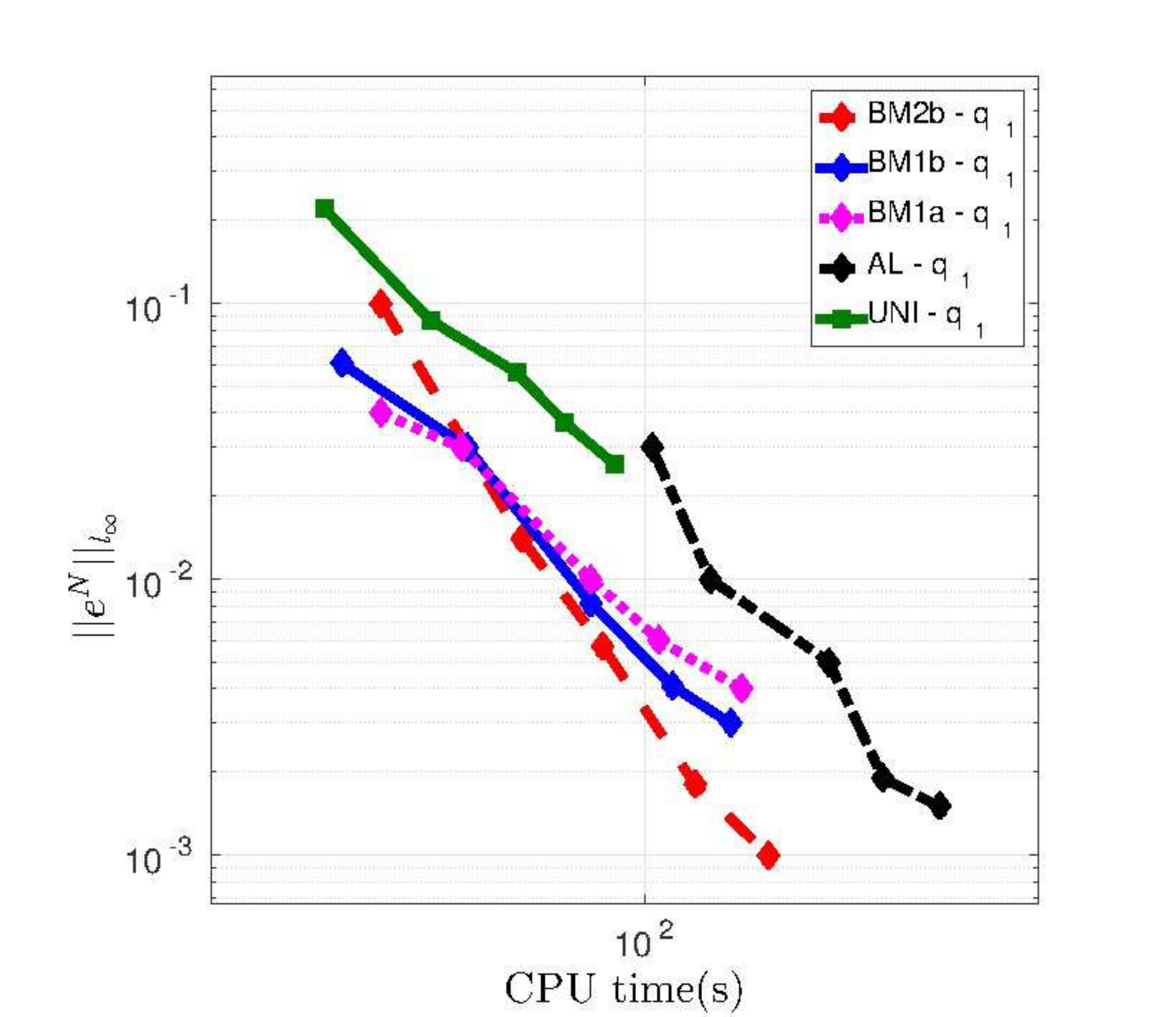}}
\scalebox{0.2}{\includegraphics{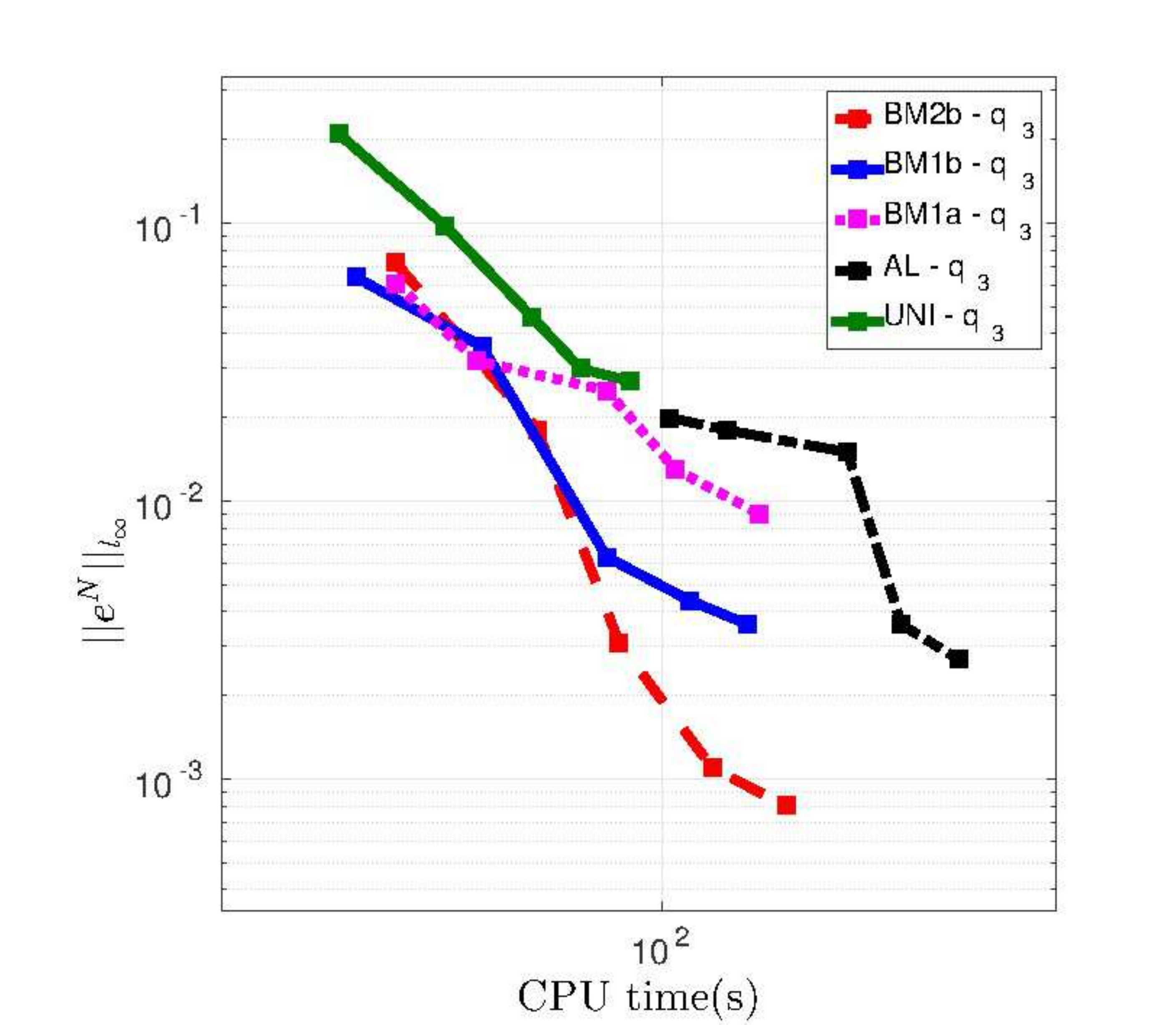}}
\scalebox{0.2}{\includegraphics{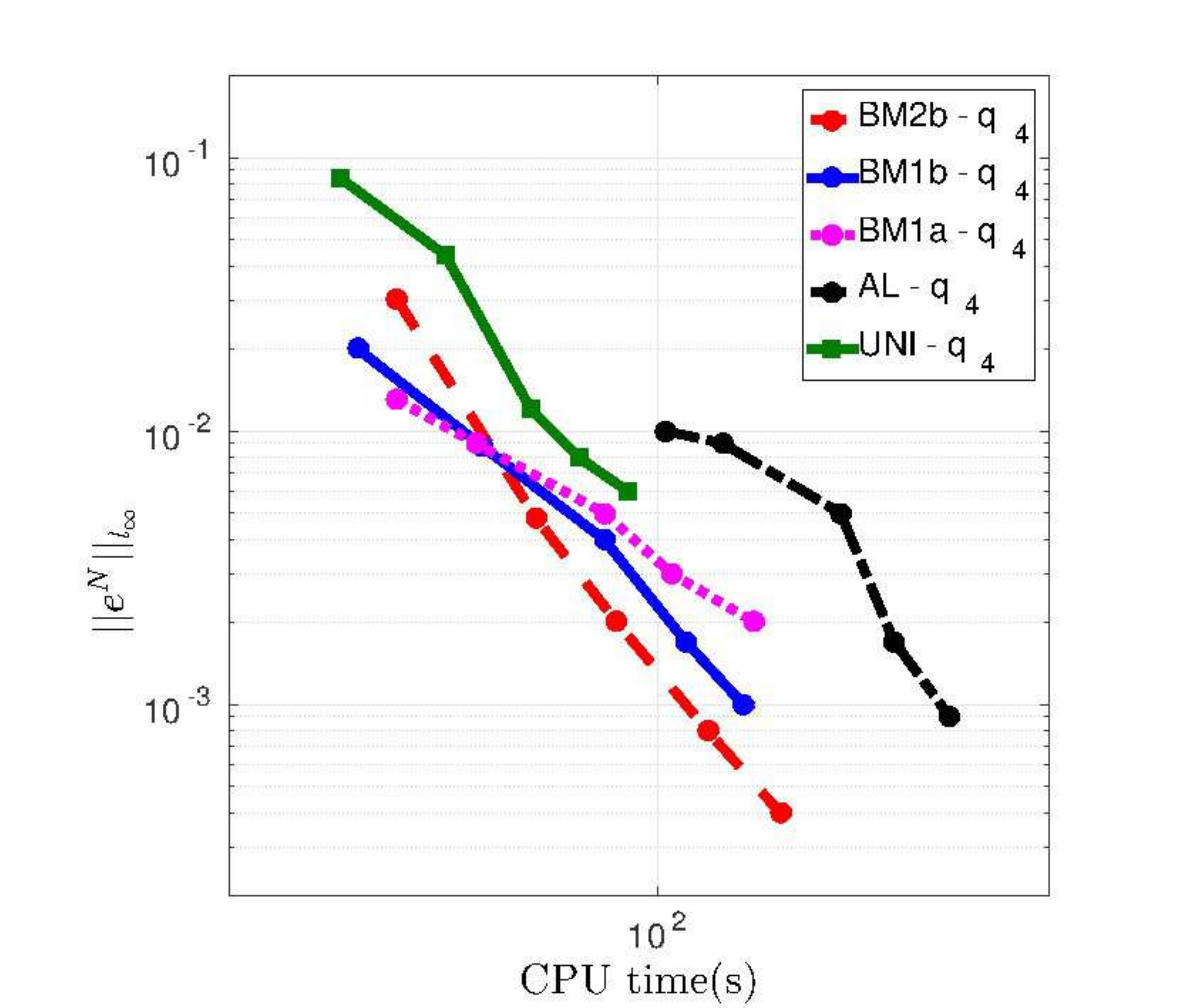}}\\
\end{minipage}
}
\caption{The $l_{\infty}$ error in $q_1,q_2,q_3$ plotted against the 
total CPU time
in seconds for each method, measured at time $t^{\ast}=0.2$ ms. The data points
correspond to grids using 122, 162, 218, 286, and 342 quadratic triangular 
elements. \label{fig:cost}}
\end{figure}

In addition to the accuracy of approximations produced, we must take into 
account the computational cost using each monitor function. The plots in Figure~\ref{fig:cost} show 
the $l_{\infty}$ error (\ref{linf}) for the three non-zero components
of the ${\bfa Q}$-tensor ($q_1$, $q_3$, $q_4$) plotted against the total CPU time in seconds 
required for each method.
Also included here for comparison are the results with uniform meshes with the same number of 
elements. The first important observation from these results is that, regardless of which monitor 
function is used, the MMPDE-based adaptive methods always outperform a standard uniform mesh in 
terms of this measure of efficiency. Furthermore, the results also show clearly that, as more 
accurate solutions are sought, the BM2b monitor function proves to be the most cost-effective 
choice, as in the cases other methods are cheaper, the error is unattractively large.
Although the more traditional arc-length based monitor function (AL) comes 
closest to matching the accuracy of BM2b, it does so at a far greater cost. 
Hence, overall, we conclude that the BM2b combination of monitor and input functions is the method 
of choice.

\subsection{Test problem 2: two-dimensional Pi-cell problem}
\label{Picell}
\begin{figure}[!ht]
\begin{center}
\includegraphics[width=12cm,height=4cm]{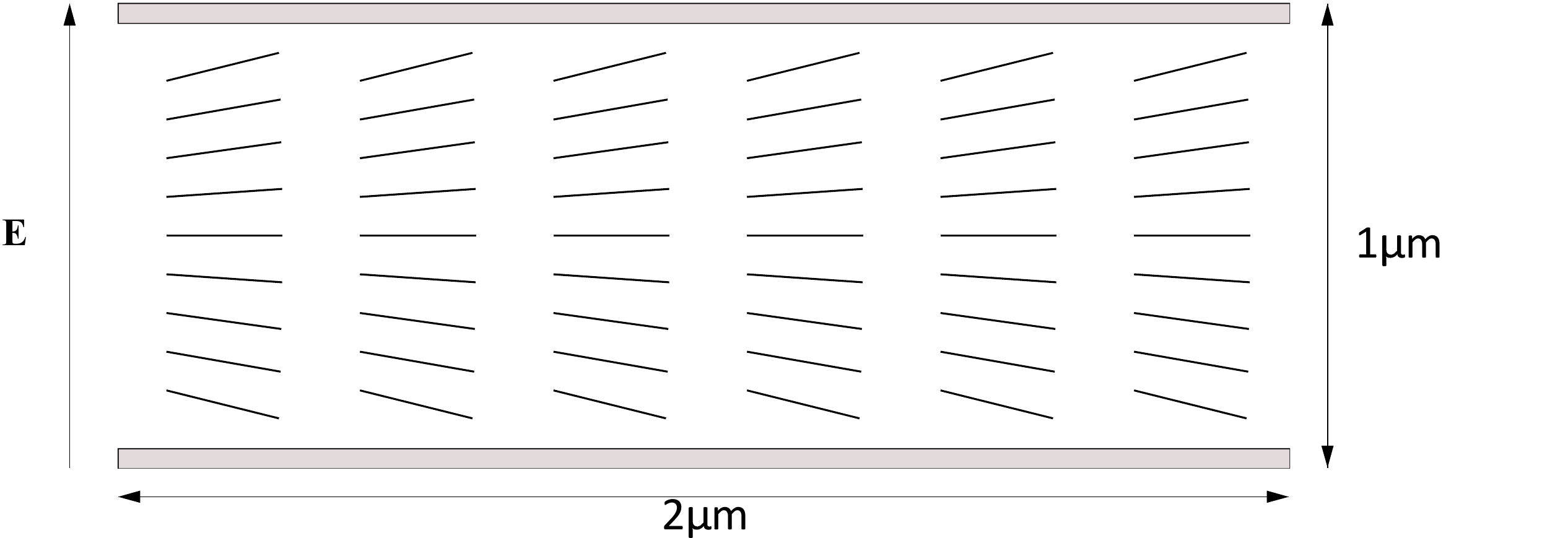}\\
(a) Schematic of unperturbed Pi-cell problem.\\
\includegraphics[width=12cm,height=4cm]{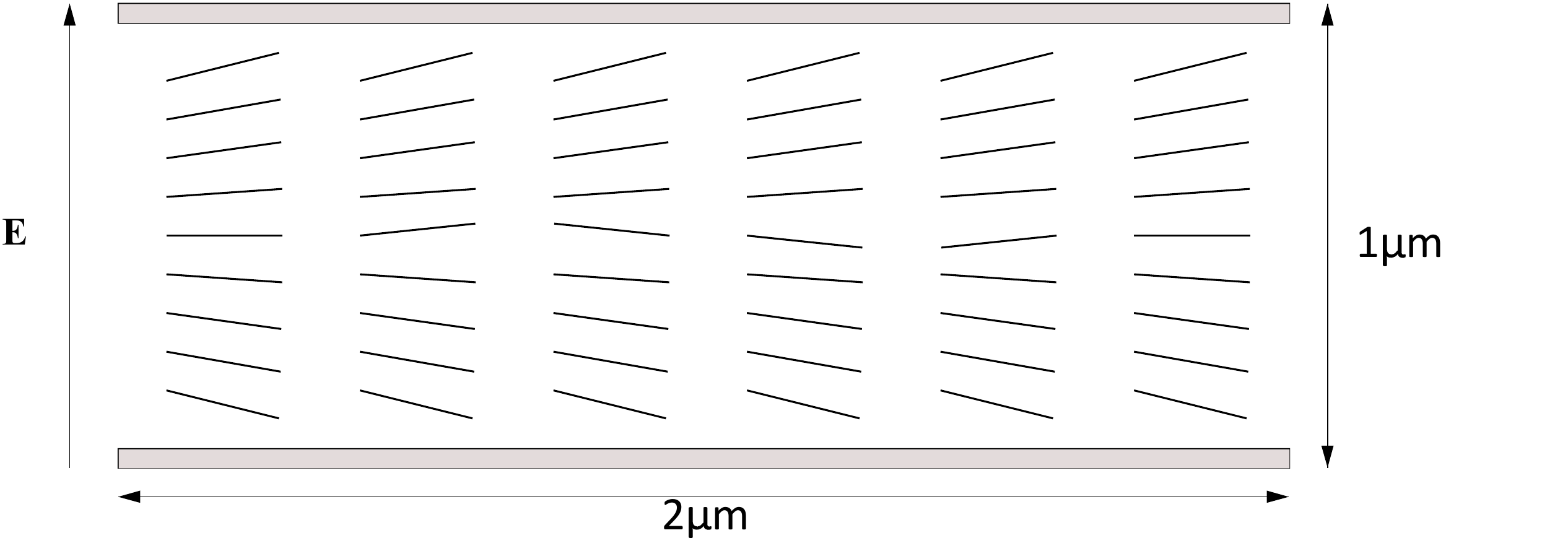}
(b) Schematic of sinusoidally perturbed Pi-cell problem.\\
\end{center}
\caption{Cell configuration in splay state \ar{under the influence 
of an electric field $\mathbf{E}$.}
\label{pipic}}
\end{figure}

We now consider a fully two-dimensional time-dependent problem, involving defects which move in 
time through a liquid crystal cell, eventually annihilating each other to leave an unperturbed 
state. 
The geometry is that of a Pi-cell \cite{boyd:80} of width two microns and thickness
one microns, with liquid crystal parameters again taken from \cite{barberi:04}.
At both boundaries, the cell surface is
treated so as to induce alignments uniformly tilted by a specified
tilt angle, $\theta_t$, but oppositely directed. If a sufficiently high voltage is applied 
across the Pi-cell for long enough, then a transition 
from the splay state (which has mostly horizontal alignment of the director with a slight splay, as 
depicted in Figure~\ref{pipic}(a))
to the bend state (which has mostly vertical alignment of the director with a bend of almost $\pi$ 
radians) can be achieved. Based on experimental results, two different physical mechanisms for this 
transition have been proposed:
a homogeneous transition via the material melting in the central plane of the Pi-cell,
or an inhomogeneous transition mediated by the nucleation of defect pairs which move and 
eventually annihilate each other. The homogeneous transition problem is essentially 
one-dimensional and has previously been modelled by several authors using moving mesh techniques 
\cite{abl10,abl11,macdonald:11,rn08}. For a more challenging test of our two-dimensional adaptive 
moving mesh approach, we will concentrate here on the simulation of the inhomogeneous transition 
type with moving defects.  This problem is still in theory relatively unchallenging: at $t=0$, if no 
perturbation is applied, the director angle simply varies linearly between the tilt angles, as in 
Figure~\ref{pipic}(a), with a director angle across the middle of the cell of $\theta= 0^\circ$. 
In practice, however, it is unrealistic for this to be achieved exactly in a physical cell due to 
small variations in the pretilt angles or thermal fluctuations.
We therefore follow \cite{bos:07} and modify the initial director angle across the middle of the 
cell so that it follows the sinusoidal function $\sin(2\pi x/p)$, where 
$x$ is the spatial coordinate in the horizontal plane, and $p$ is the cell width.
This perturbation is fixed only at $t=0$ for one time step,
but introduces solution gradients in two dimensions, as portrayed in Figure \ref{pipic}(b), which 
provide a bigger challenge for our numerical method. 

We consider a cell of width 2$\mu$m 
and thickness 1$\mu$m, with a pre-tilt angle of $\theta=\pm 6^\circ$. 
An electric field of strength $18V{\mu}{\rm m}^{-1}$ is applied
parallel to the cell thickness at time $t=0$.
Based on the evidence from the static defect problem in \S\ref{static}, we present results 
from the BM2b monitor function only.

\begin{figure}[!ht]
\begin{center}
\hspace*{-1.9cm}
\includegraphics[width=16cm,height=4cm]{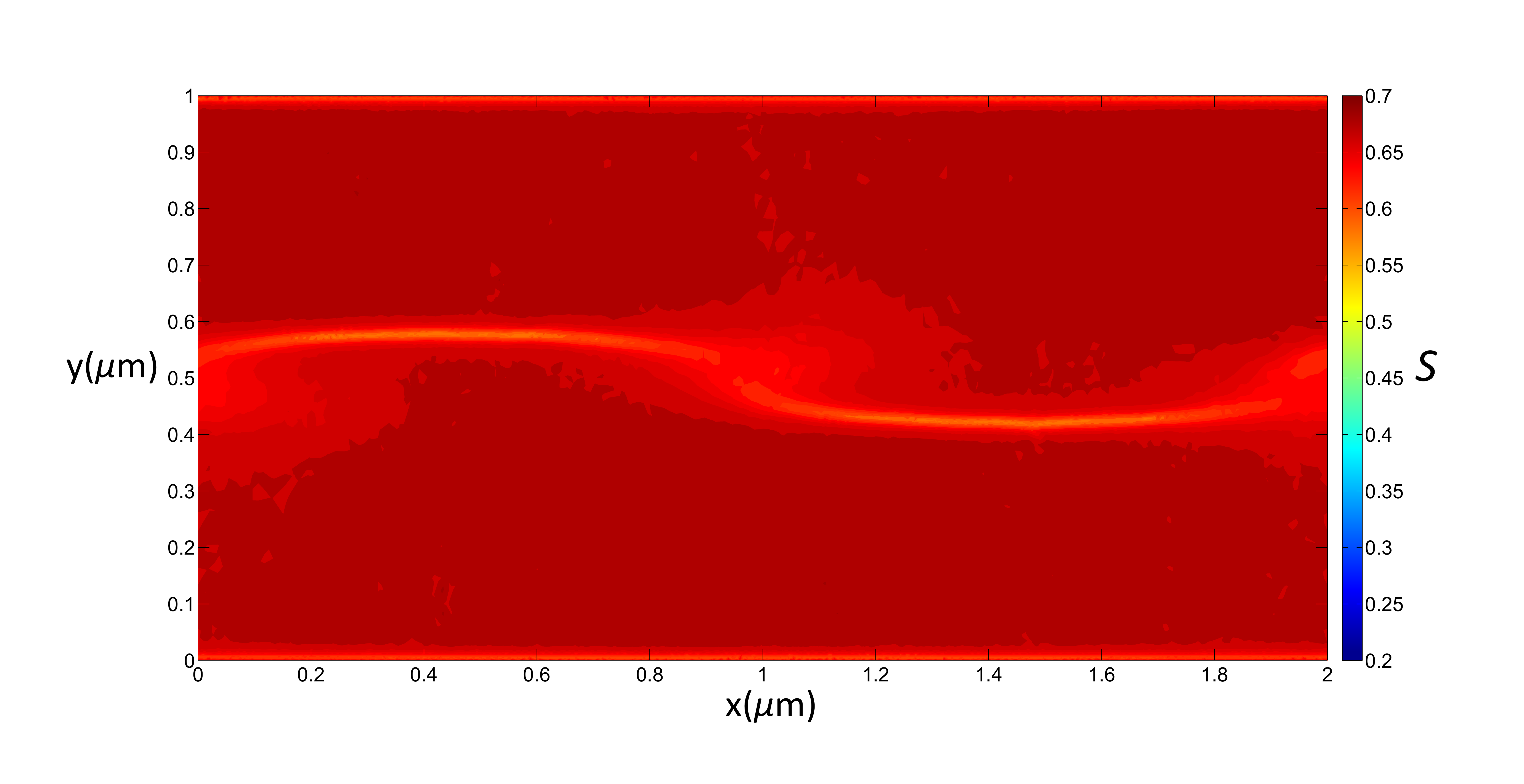}\\
(a) {Order parameter.}\\
\hspace*{-1.9cm}
\includegraphics[width=16cm,height=4cm]{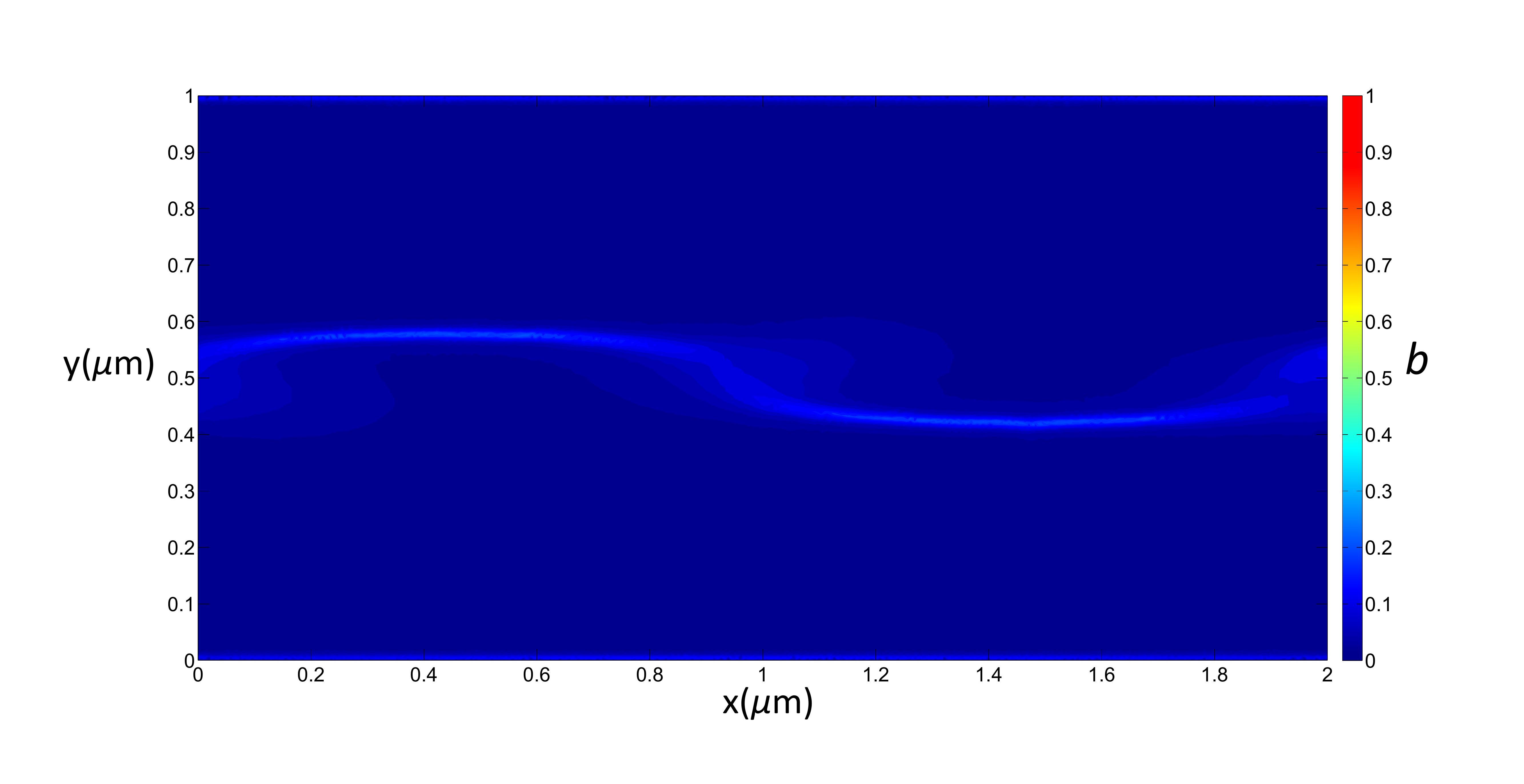}\\
(b) {Biaxiality.}\\
\hspace*{-1.9cm}
\includegraphics[width=16cm,height=4cm]{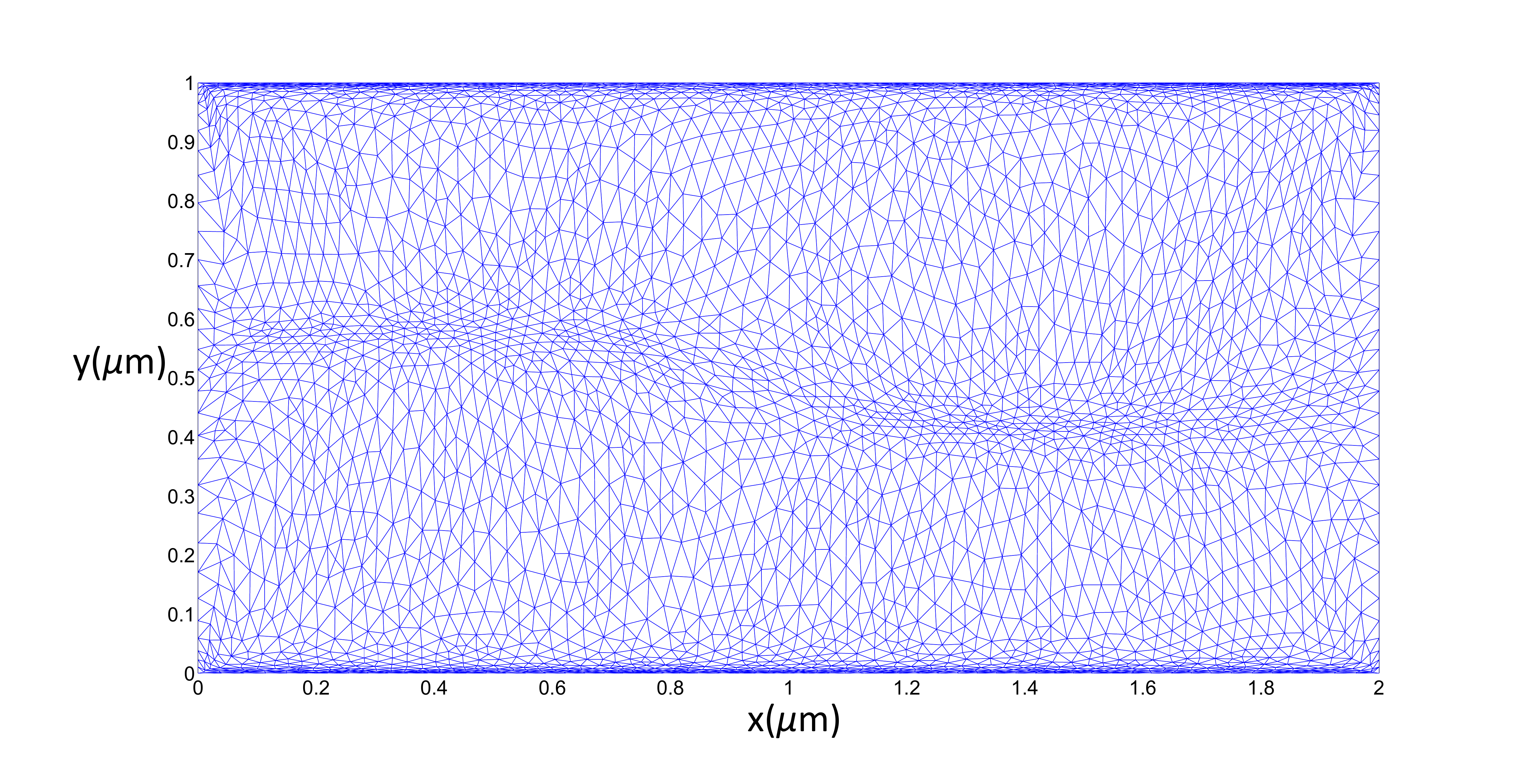}\\
(c) {Mesh.}
\end{center}
\caption{Order parameter, biaxiality and adapted mesh after 12$\mu$s.
\label{pi_ord_bi}}
\end{figure}

\begin{figure}[!htbp]
\centering
\hspace*{-1.2cm}
\includegraphics[width=13.6cm,height=5.4cm]{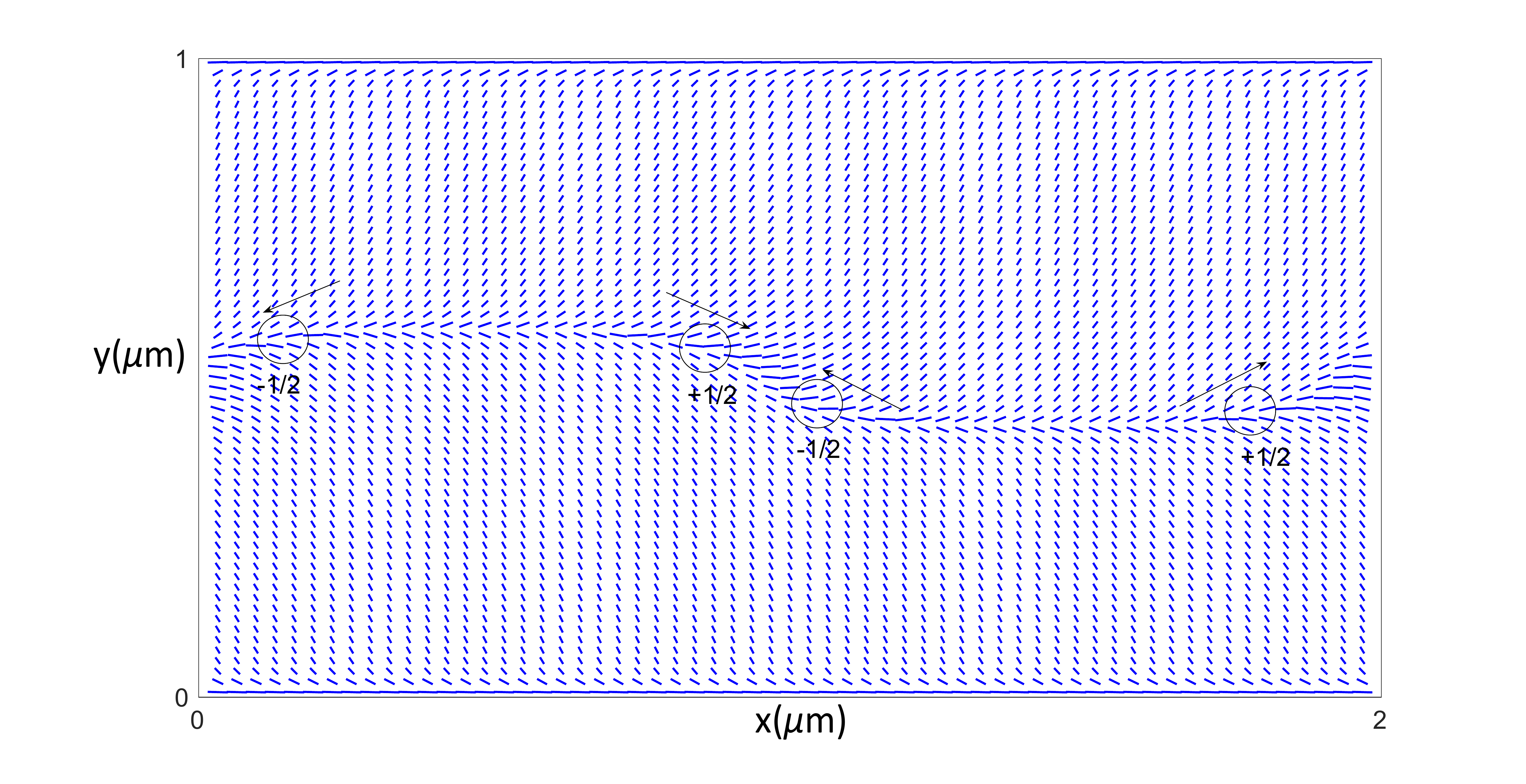}
\caption{Director field after 15.5$\mu$s. \label{bos_15d}}
\end{figure}

\begin{figure}[!htbp]
\centering
\mbox{
\begin{minipage}{5.2in}
\hspace{-1.5cm}
\scalebox{0.3}{\includegraphics{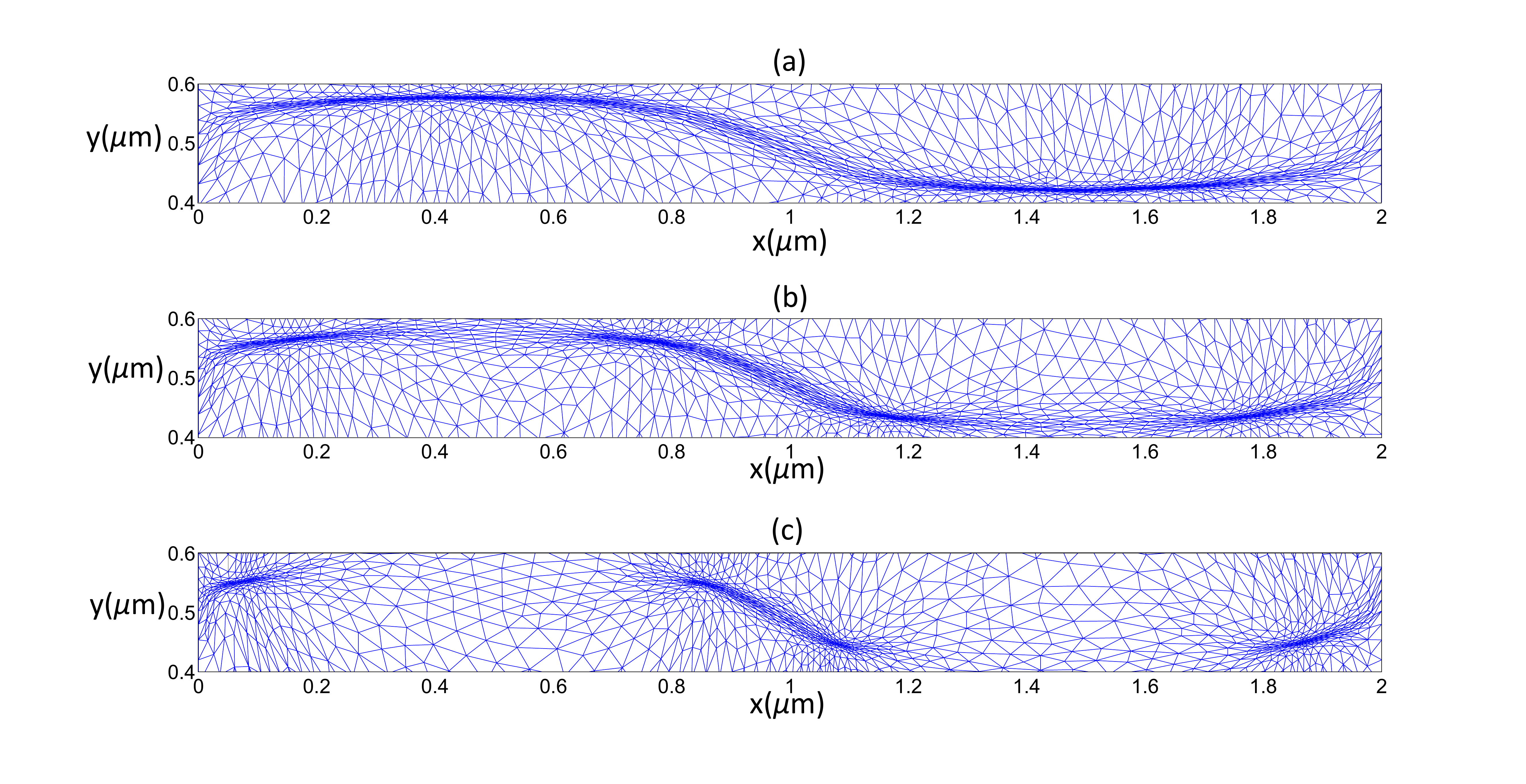}}
\caption{Adapted finite element mesh after (a) 15.5$\mu$s (b) 16$\mu$s and 
(c) 17$\mu$s.}
\label{m15_d2}
\end{minipage}}
\end{figure}

\begin{figure}[!htbp]
\centering
\mbox{
\begin{minipage}{5.2in}
\hspace{-1.5cm}
\scalebox{0.3}{\includegraphics{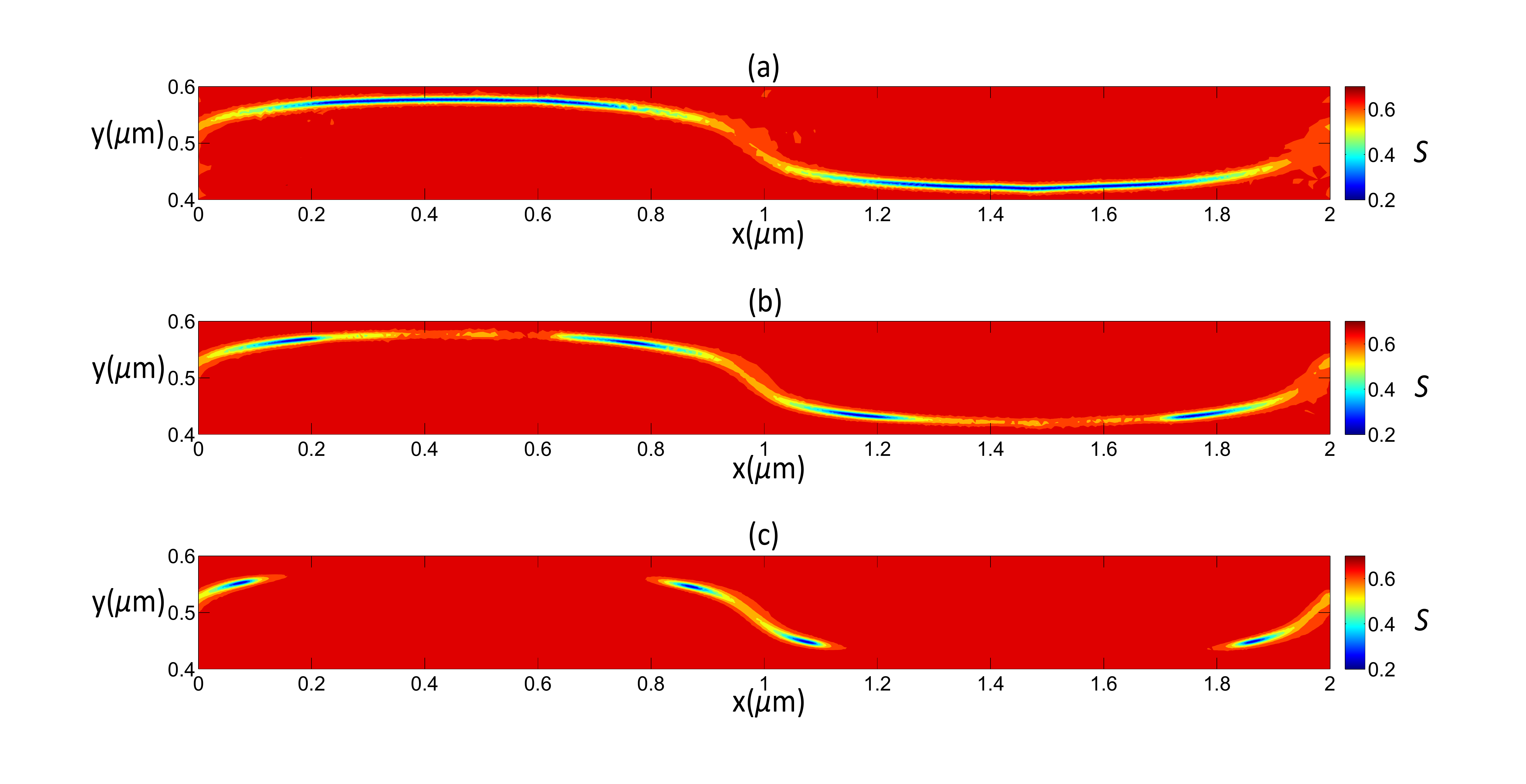}}
\caption{Order parameter profile after (a) 15.5$\mu$s (b) 16$\mu$s and 
(c) 17$\mu$s}
\label{s15_d2}
\end{minipage}}
\end{figure}

\begin{figure}[!htbp]
\centering
\mbox{
\begin{minipage}{5.2in}
\hspace{-1.5cm}
\scalebox{0.3}{\includegraphics{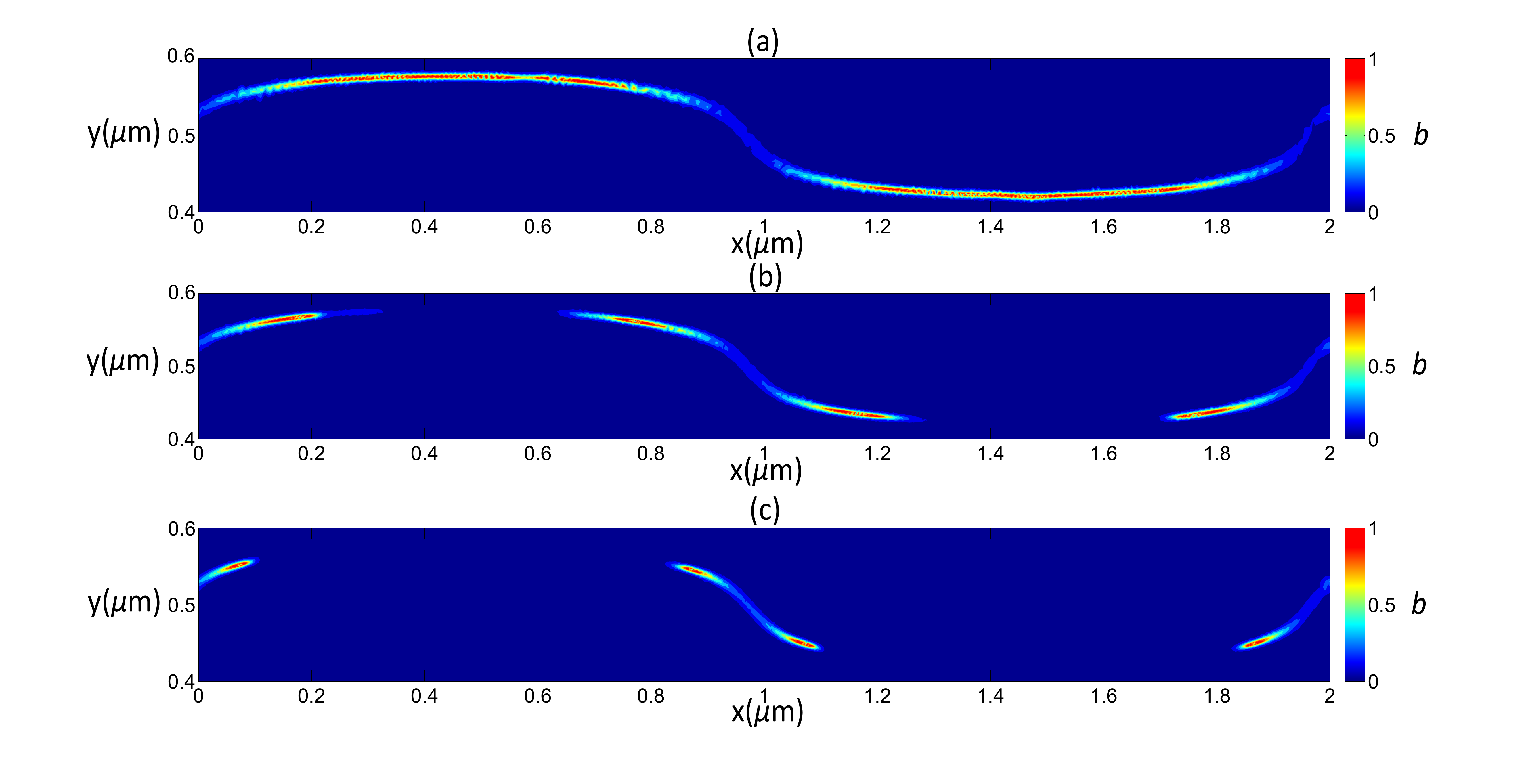}}
\caption{Biaxiality profile after (a) 15.5$\mu$s (b) 16$\mu$s and (c) 17$\mu$s}
\label{b15_d2}
\end{minipage}}
\end{figure}

Initially, immediately before the application of the electric field, the cell is 
in an equilibrium state where the order parameter and biaxiality take 
constant values of $0.65$ and $0$, respectively.  The mesh at this 
stage is quasi-uniform as no adaptivity has yet taken place. As time evolves, the 
combined effect of the perturbation and the applied electric field become 
apparent.
Figure~\ref{pi_ord_bi} shows the cell state 12$\mu$s after the 
application of the electric field; at this time there is a region of 
concentrated splay distortion at the centre of the cell.  
Within this area, the order parameter and biaxiality are no 
longer at their constant equilibrium values: the mesh, as expected, has started to adapt
as depicted in Figure \ref{pi_ord_bi}(c).
After 15.5$\mu$s the distortion at the centre of the cell has become more pronounced, and we can 
clearly observe pairs of +1/2 and -1/2 defects within this area, as shown in Figure \ref{bos_15d}.
Outwith the distorted area, the cell is largely in an equilibrium state, with the order parameter 
and biaxiality still at their constant equilibrium values. However, the cores of the defects are now 
completely biaxial, and the measure of biaxiality approaches its maximal value of 1.
Recalling that the BM2b monitor function uses biaxiality as its input, it is not surprising 
that the mesh has now adapted significantly from its quasi-uniform initial state, and has
started to adapt well to resolve the defects (see Figure \ref{m15_d2}).
As time evolves further, the oppositely
signed defects are attracted to each other, moving ever closer until they ultimately
meet and annihilate each other. Figures~\ref{s15_d2} and \ref{b15_d2} show snapshots of the order
parameter and biaxiality respectively, measured after 15.5$\mu$s, 16$\mu$s and 17$\mu$s, 
calculated on the meshes shown in Figure~\ref{m15_d2}.
In Figure \ref{s15_d2}, after 16$\mu$s, the mesh is still well adapted to the sinusodial
shape of the initial perturbation, consistent with the presence of large variations in
the biaxiality throughout the central area of the cell. However, after 17$\mu$s, by which
point the defects have almost coalesced, the mesh has relaxed in areas where the
biaxiality is now back to its equilibrium value, and instead is completely focused
on resolving the defects. After the defects meet and annihilate, the biaxiality and
order parameter again relax towards their equilibrium value everywhere in the cell,
and the mesh also relaxes back to a quasi-uniform state.
Overall, throughout the simulation, the adaptive moving mesh
method does an excellent job of tracking the development, movement, and
annihilation of the defects in the liquid crystal cell. In particular, the method
is able to cope well with the small-scale structure of the defect core, and the
short timescales associated with the establishment and annihilation of defects.

\section{Summary}

The focus of this paper is on the description and application of a new efficient
moving mesh method for $\bfa Q$-tensor models of liquid crystal cells. 
Although some of the ideas contained here are described in a one-dimensional setting in 
\cite{macdonald:15}, extending the method to tackle the more physically realistic fully 
two-dimensional problems presented here required us to address a number of significant new 
challenges. As with all moving mesh methods, the choice of an appropriate adaptivity criterion is 
crucial: here we have established that using a monitor function based on second-order partial 
derivatives of a local measure of the biaxiality of the liquid crystal material (BM2b in 
Table{\ref{montab}) is extremely successful in this regard. Using a test problem based on a 
static +1/2 defect, we demonstrated in \S\ref{static} computed solutions from all of our 
proposed methods converge optimally in space. However, a comparison of efficiency demonstrated that 
the BM2b monitor function is clearly the method of choice in terms of computationally efficiency 
when a reasonable level of accuracy is required. Furthermore, when applied to the more realistic 
but more challenging fully two-dimensional Pi-cell problem described in \S\ref{Picell},
the adaptive MMPDE method based on the BM2b monitor function proved to be very effective for 
resolving the movement and core details of defects, including the
creation and annihilation of these moving singularities. This is particularly impressive given the 
very short length and time scales involved in these aspects of the material's behaviour.

Of course, some challenges still remain. Particularly useful in practice would be the extension of 
our method to multi-dimensional problems with irregular geometries. This 
would pose a further challenge to the adaptive moving mesh method
as it would potentially have to resolve defects present around the areas where the
cell geometry is most complex. A prime candidate would be the Zenith Bistable-Device (ZBD) 
described in \cite{bryanbrown:95,newton:97}, where the liquid crystal 
cell has an alignment layer on the upper surface and a periodic grating structure on the
lower surface. 

\section*{Acknowledgements}
The authors would like to thank Dr Chris Newton, formerly of Hewlett-Packard 
plc, for invaluable discussions. This work was supported in part by an Engineering and Physical 
Sciences Research Council Industrial CASE Studentship (EP/P505747/1) with Hewlett-Packard plc.\\

\bibliography{MMLC}

\end{document}